\documentclass{amsart}

\usepackage{hyperref}
\usepackage{amsmath,amsthm,amssymb}
\theoremstyle{definition}
\newtheorem{theorem}{Theorem}
\newtheorem{proposition}{Proposition}
\usepackage{mathrsfs}
\usepackage{amsfonts}
\usepackage{graphicx}
\usepackage{epstopdf}
\usepackage{algorithmic}
\usepackage{subfigure}
\usepackage{cleveref}
\allowdisplaybreaks[1]

\ifpdf
\hypersetup{
  pdftitle={Construction of conformal maps based on the locations of singularities
	for improving the double exponential formula},
  pdfauthor={Shunki Kyoya and Ken'ichiro Tanaka}
}
\fi

\begin{document}
\title
[Integration with conformal maps]
{Construction of conformal maps based on the locations of singularities
for improving the double exponential formula}
\author{Shunki Kyoya}
\address{Department of Mathematical Informatics,
	Graduate School of Information Science and Technology,
	The University of Tokyo,
	7-3-1, Hongo, Bunkyo-ku, Tokyo 113-8656, Japan} 
\email{kyoya\_shunki@mist.i.u-tokyo.ac.jp}
\author{Ken'ichiro Tanaka}
\address{Department of Mathematical Informatics,
	Graduate School of Information Science and Technology,
	The University of Tokyo,
	7-3-1, Hongo, Bunkyo-ku, Tokyo 113-8656, Japan} 
\email{kenichiro@mist.i.u-tokyo.ac.jp}

\begin{abstract}
	The double exponential formula, or the DE formula,
	is a high-precision integration formula
	using a change of variables called a DE transformation;
	whereas there is a disadvantage 
	that it is sensitive to singularities of an integrand near the real axis.
	To overcome this disadvantage,
	Slevinsky and Olver (SIAM J. Sci. Comput., 2015) attempted to improve it
	by constructing conformal maps
	based on the locations of singularities.
	Based on their ideas, we construct a new transformation formula.
	Our method employs special types of the Schwarz-Christoffel transformations
	for which we can derive their explicit form. 
	Then, the new transformation formula
	can be regarded as a generalization of the DE transformations.
	We confirm its effectiveness by numerical experiments.
\end{abstract}

\maketitle

%%%%%%%%%%%%%%%%%%%%%%%%%%%%%%%%%%%%%%%%%%%%%%%%%%%%%%%%%%%%%%%%%%%%%%%%%%%
\section{Introduction}
\label{sec:1}
%%%%%%%%%%%%%%%%%%%%%%%%%%%%%%%%%%%%%%%%%%%%%%%%%%%%%%%%%%%%%%%%%%%%%%%%%%%

\textit{The double exponential formula}, or \textit{the DE formula},
is a numerical integration formula 
using a change of variables called a DE transformation
and the trapezoidal rule \cite{DE}.
For example, an integral on the interval $(-1, 1)$ is calculated
as follows:
\begin{align}
	\int_{-1}^1 f(x) \mathrm{d} x \approx
	h \sum_{j=-n}^n f(\phi(jh))\phi'(jh),
\end{align}
where the DE transformation corresponding to this interval is
\begin{align}
	\phi(t) = \tanh \left( \frac{\pi}{2} \sinh(t) \right).
\end{align}
\textit{DE transformations} are changes of variables
which make the transformed integrands decay
double exponentially:
\begin{align}
	f(\phi(t))\phi'(t) = \mathrm{O}(\exp(-\beta\mathrm{e}^{|t|}))
	\quad (t \to \pm\infty)
\end{align}
for some $\beta>0$.

The advantages and disadvantages of the DE formula have been formulated
by Sugihara \cite{Sugihara}.
He estimated them using a parameter $d$,
the width of the domain around the real axis in which the transformed
integrand is analytic.
He has shown that the error of the DE formula converges on the order of
$\mathrm{O}(\mathrm{e}^{-kN/\log N})$ as $N \to \infty$,
where $N$ is the number of the nodes for the trapezoidal rule
and $k$ is proportional to $d$.
From this formulation, we see that the DE formula makes the error converges rapidly
regardless of end-point singularities;
whereas it has a disadvantage that is sensitive to singularities of the
integrand near the real axis
since they make the parameter $d$ small.

To overcome this problem,
Slevinsky and Olver \cite{SO} proposed to improve the DE formula 
by modifying the DE transformations.
They have revealed relations between the parameter $d$ and singularities,
and proposed to make polynomial adjustments to the DE transformations
based on the locations of singularities.

In this paper, we construct a new transformation formula
based on their idea.
First, we list the options of the transformation formulas
using the idea of the Schwarz-Christoffel transformation.
Then, we choose the optimal one
from the perspective of the precision and ease of
numerical integration.
This transformation formula is not only a modification of the DE transformations,
but also can be considered to be a generalization of them.
We confirm its effectiveness by numerical experiments.

The rest of the paper is organized as follows.
In \Cref{sec:2}, we summarize the Sugihara's analysis.
In \Cref{sec:3}, we describe the idea to improve the DE formula
and introduce the method by Slevinsky and Olver.
In \Cref{sec:4}, we show the proposed methods.
In \Cref{sec:5}, we show numerical experiments.
Finally, we conclude this paper in \Cref{sec:6}.

We show proofs and calculations omitted in this paper in the appendix.
Programs for the proposed methods are available from \cite{geDEprogram}.

%%%%%%%%%%%%%%%%%%%%%%%%%%%%%%%%%%%%%%%%%%%%%%%%%%%%%%%%%%%%%%%%%%%%%%%%%%%
\section{Precision of  the double exponential formula}
\label{sec:2}
%%%%%%%%%%%%%%%%%%%%%%%%%%%%%%%%%%%%%%%%%%%%%%%%%%%%%%%%%%%%%%%%%%%%%%%%%%%

On the basis of theorems in \cite{Sugihara},
we formulate the precision of the DE formula
by evaluating the error of the trapezoidal formula
in the case where the integrands decay double exponentially.

We define a family of integrands. Let $d$ be a positive number and 
let $\mathscr{D}_d$ denote the strip region of width $2d$ 
about the real axis:
\begin{align}
	\mathscr{D}_d = \{ z \in \mathbb{C} \mid |\mathrm{Im} z| < d \}.
\end{align}
Let $\omega$ be a non-vanishing function defined on the region $\mathscr{D}_d$, 
and define the \textit{weighted Hardy space} $H^{\infty}(\mathscr{D}_d, \omega)$ by
\begin{align}
	H^{\infty}(\mathscr{D}_d, \omega)
	= \{ f: \mathscr{D}_d \to \mathbb{C} \mid 
	f(z)\ \text{is analytic in}\ \mathscr{D}_d,\ \text{and}\ ||f|| < \infty \},
\end{align}
where the norm of $f$ is given by
\begin{align}
	||f|| = \sup_{z \in \mathscr{D}_d}
	\left| f(z) / \omega(z) \right|.
	\label{eq:normf}
\end{align}

For the following discussions, we assume that the function $\omega$
decays double exponentially.
Then, since
\begin{align}
	|f(z)| \leq ||f||\,|\omega(z)| \quad (z \in \mathscr{D}_d)
\end{align}
holds from the definition of the norm \cref{eq:normf},
the weighted Hardy space $H^{\infty}(\mathscr{D}_d, \omega)$
represents the family of integrands which decay double exponentially.

Let $N = 2n+1$ be the number of the nodes for numerical integration.
The error of the $N$-point trapezoidal formula
is estimated using an error norm.
Let $\mathscr{E}_{N,h}^{\mathrm{T}} (H^\infty (\mathscr{D}_d, \omega))$
denote the error norm in $H^\infty (\mathscr{D}_d, \omega)$:
\begin{align}
	\mathscr{E}_{N,h}^{\mathrm{T}} (H^\infty (\mathscr{D}_d, \omega)) = 
	\sup_{\substack{f \in H^\infty (\mathscr{D}_d,\,\omega) \\ ||f|| \leq 1}}
	\left|\int_{-\infty}^\infty f(x)\mathrm{d}x - h\sum_{j=-n}^n f(jh)
	\right| .
\end{align}

The following theorem gives the upper bound of this error norm.
Let $B(\mathscr{D}_d)$ denote the family of functions $g$
which are analytic in $\mathscr{D}_d$ and satisfy
\begin{align}
	\int_{-d}^d |g(x+\mathrm{i}y)|\mathrm{d}y \to 0 \quad (x \to \pm \infty)
\end{align}
and
\begin{align}
	\lim_{y \to d-0} &\int_{-\infty}^\infty
	(|g(x+\mathrm{i}y)| + |g(x-\mathrm{i}y)| ) \mathrm{d}x < \infty.
\end{align} 
\begin{theorem}[Sugihara \cite{Sugihara}] \label{thm:Sugihara_DE}
	Suppose that the function $\omega$
	satisfies the following three conditions:
	\begin{enumerate}
		\item
		$\omega \in B(\mathscr{D}_d)$;
		\item
		$\omega$ does not vanish at any point in $\mathscr{D}_d$
		and takes real values on the real axis;
		\item
		the decay rate on the real axis of $\omega$ is specified by 
		\begin{align} \nonumber
			\alpha_1\exp (-\beta_1 \mathrm{e}^{\gamma |t|}) \leq
			|\omega(t)| \leq
			\alpha_2\exp (-\beta_2 \mathrm{e}^{\gamma |t|}), 
		\end{align}
		where $\alpha_1$, $\alpha_2$, $\beta_1$, $\beta_2$, $\gamma>0$.
	\end{enumerate}
	Then the upper bound of the error norm is given as
	\begin{align}
		\mathscr{E}_{N,\,h}^{\mathrm{T}} (H^\infty (\mathscr{D}_d, \omega))
		\leq
		C_{d,\,\omega} \exp\!
			\left(-\frac{\pi d \gamma N}{\log(\pi d \gamma N / \beta_2)}\right), 
		\label{eq:error}
	\end{align}
	where $N=2n+1$, 
	$C_{d, \omega}$ is a constant
	depending on $d$ and $\omega$, 
	and the mesh size $h$ is chosen as
	\begin{align}
		h = \frac{\log(2\pi d \gamma n / \beta_2)}{\gamma n}. \label{eq:hDE}
	\end{align}
\end{theorem}
This theorem shows
that the error of the trapezoidal formula converges exponentially
according the parameters $d, \gamma$ and $\beta_2$.
The larger these parameters are, the better the convergence rate becomes.
However, it has also been shown that the parameters have a restriction as shown by
the following theorem.
\begin{theorem}[Sugihara \cite{Sugihara}] \label{thm:Sugihara_nexists}
	There exists no function $\omega$
	that satisfies the following three conditions:
	\begin{enumerate}
		\item
		$\omega \in B(\mathscr{D}_d)$;
		\item
		$\omega$ does not vanith at any point in $\mathscr{D}_d$
		and takes real values on the real axis;
		\item
		the decay rate on the real axis of $\omega$ is specified by
		\begin{align} \nonumber
			\omega(t) = \mathrm{O}(\exp(-\beta\exp(\gamma|t|)))
			\quad as\ |t| \to \infty, t\in \mathbb{R},
		\end{align}
		where $\beta>0$ and $\gamma>\pi/(2d)$.
	\end{enumerate}
\end{theorem}

DE transformations are changes of variables 
which make the transformed integrand $f(\phi(\cdot))\phi'(\cdot)$
be a member of $H^\infty (\mathscr{D}_d, \omega)$.
However, they do not necessarily guarantee the optimal
convergence rate of \Cref{thm:Sugihara_DE}.
In the following sections, we improve the DE formula
by modifying the DE transformations
so that the convergence rate will be better.

%%%%%%%%%%%%%%%%%%%%%%%%%%%%%%%%%%%%%%%%%%%%%%%%%%%%%%%%%%%%%%%%%%%%%%%%%%%
\section{Improvement of the double exponential formula}
\label{sec:3}
%%%%%%%%%%%%%%%%%%%%%%%%%%%%%%%%%%%%%%%%%%%%%%%%%%%%%%%%%%%%%%%%%%%%%%%%%%%

We consider DE transformations which are written as
$\phi(t) = \psi\!\left(\frac{\pi}{2}\sinh(t) \right)$ for some function $\psi$.
We show examples of such DE transformations in \Cref{tab:DE},
where the fourth was introduced in \cite{DEindifinite}.
In these cases,
$\psi$ are periodic with period $2\pi$ in the direction of the imaginary axis.
We improve the DE formula by constructing a function $H$
and changing these transformations to $\phi(t) = \psi(H(t))$.

On the basis of theorems in \Cref{sec:2}, Slevinsky and Olver \cite{SO}
proposed an idea to construct $H$ appropriately.
In this section, we describe their idea and method from our point of view.

\begin{table}[th]
{\footnotesize
	\caption{Examples of DE transformations
	which are written as $\phi(t) = \psi\!\left(\frac{\pi}{2}\sinh(t) \right)$.
	\label{tab:DE}}
	\begin{center}
		\begin{tabular}{cccc} \hline
			interval & integrand & $\phi(t)$ & $\psi$ \\ \hline
			$(-1, 1)$ & $f(x)$ & $\tanh(\frac{\pi}{2}\sinh(t))$ & $\tanh(\cdot)$ \\
			$(-\infty, \infty)$ & $f(x)$ & $\sinh(\frac{\pi}{2}\sinh(t))$ & $\sinh(\cdot)$ \\
			$(0, \infty)$ & $f(x)$ & $\exp(\frac{\pi}{2}\sinh(t))$ & $\exp(\cdot)$ \\
			$(0, \infty)$ & $f_1(x)\mathrm{e}^{-vx}\ (v > 0)$
			& $\log(\exp(\frac{\pi}{2}\sinh(t)+1))$  & $\log(\exp(\cdot + 1))$\\ \hline
		\end{tabular}
	\end{center}
}
\end{table}

%%%%%%%%%%
\subsection{Construction of the transformation formula}
\label{subsec:H}
%%%%%%%%%%

We construct the transformation formula $H$ so that
the transformed integrand $f(\psi(H(\cdot)))\psi'(H(\cdot))H'(\cdot)$
will be a member of $H^{\infty}(\mathscr{D}_d, \omega)$.
Then, the convergence rate of \Cref{thm:Sugihara_DE} is applied.
We wish to maximize it under the limit of \Cref{thm:Sugihara_nexists}.
Thus, we wish to determine the function $H$
according to the following optimization problem:
\begin{align} \label{eq:maximize}
\max_{H}&\quad \frac{\pi d \gamma N}{\log(\pi d \gamma N/\beta_2)}
\ (\text{convergence rate of \Cref{thm:Sugihara_DE}}) \\
\text{subject to} &\quad \left\{
	\begin{aligned}
		& f(\psi(H(\cdot)))\psi'(H(\cdot))H'(\cdot) \in H^{\infty}(\mathscr{D}_d, \omega) \\
		& d > 0 \\
		& \omega\ \text{satisfies the conditions of \Cref{thm:Sugihara_DE}} \\
		& d\gamma \leq \pi/2 \ (\text{limitation of \Cref{thm:Sugihara_nexists}}).
	\end{aligned}\right.
\end{align}
However, it is difficult to solve this optimization problem generally.
In order to make the problem more simply,
we consider the asymptotic form of \cref{eq:maximize}
as $N \to \infty$.
Since it is written asymptotically as
\begin{align}
	\frac{\pi d \gamma N}{\log(\pi d \gamma N/\beta_2)}
	\approx \frac{\pi d\gamma N}{\log N} \quad (N \to \infty),
	\label{eq:error_asympt}
\end{align}
the value of $d\gamma$ is dominant.
Thus, we construct the function $H$ according to the following method:
\begin{itemize}
	\item
	First, we restrict the options of $H$ so that $d\gamma = \pi/2$ will be satisfied.
	Here, we assume that $d = \pi/2$ and $\gamma = 1$.
	\item
	Then, we choose $H$ from these options so that
	the parameter $\beta_2$ will be larger.
\end{itemize}

\begin{figure}[ht]
	\centering
	\includegraphics[scale = 0.3]{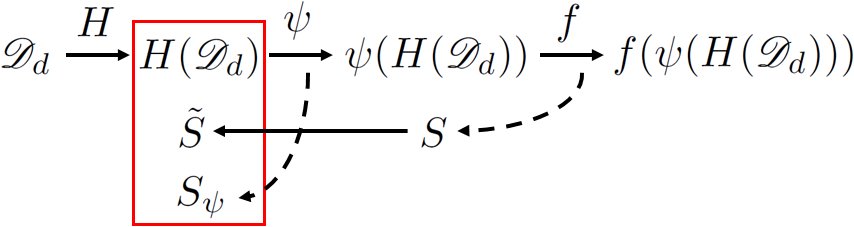}
	\caption{Relations between $H$ and singularities which $H$
	needs to avoid \label{img:sing}}
\end{figure}

The condition $d = \pi/2$ is equivalent to the condition that
the transformed integrand is analytic in $\mathscr{D}_{\pi/2}$.
It is attained by avoiding singularities.
We assume that $f$ has a finite number of 
singularities which are symmetric with respect to the real axis.
We write these singularities as
$S = \{ \delta_j \pm \epsilon_j \mathrm{i} \mid j = 1 \dots m \}$.
Let $\tilde{S}$ denote the preimage of $S$ by $\psi$.
By the periodicity of $\psi$, we write elements of $\tilde{S}$ as
\begin{align}
	\tilde{S} = \{ \tilde{\delta}_j \pm
	(\tilde{\epsilon}_j + 2k\pi) \mathrm{i}
	\mid j = 1, \ldots, m,\, k \in \mathbb{Z} \},
\end{align}
where $\tilde{\delta}_1 < \dots < \tilde{\delta}_m$
and $0 < \tilde{\epsilon}_j \leq \pi\,(j = 1, \ldots, m)$.
Also, if  $\psi$ has singularities, we write them as $S_\psi$.
For example, $\psi = \tanh(\cdot)$ has singularities
$S_\psi = \{ ( \pm \pi/2 + 2k\pi)\mathrm{i} \mid k \in \mathbb{Z} \}$.
In order to make the transformed integrand
analytic in $\mathscr{D}_{\pi/2}$,
we need to make the image $H(\mathscr{D}_{\pi/2})$
avoid the singularities in $\tilde{S}$ and $S_\psi$;
that is, we construct the function $H$
so that it will satisfy
\begin{align} \label{eq:sing}
		s \notin H (\mathscr{D}_{\pi/2}) \quad (s \in \tilde{S} \cup S_\psi).
\end{align}
\Cref{img:sing} summarizes this condition.

The parameters $\gamma$ and $\beta_2$ appears in the coefficients
of the transformation formula $H$.
For example, if it is written as
$H(t) = C\sinh(\gamma't) + \mathrm{o}(\mathrm{e}^{\gamma't})$ as $|t| \to \infty$,
then we see that $\gamma' = \gamma$ and
that $C$ is proportional to $\beta_2$.
Thus, we fix $\gamma'$ to 1 and consider
how to make $C$ larger under the condition \cref{eq:sing}.

%%%%%%%%%%
\subsection{DE transformations}
\label{subsec:HDE}
%%%%%%%%%%

We review the DE transformations
in the context of \Cref{subsec:H}.
We can write the function $H$ of the DE transformations by
\begin{align}
	H_{\mathrm{DE}}(t) = \frac{\pi}{2}\sinh(t).
\end{align}
The image $H_{\mathrm{DE}}(\mathcal{D}_{\pi/2})$ is
the entire complex plane with a pair of slits:
\begin{align}
	H_{\mathrm{DE}}(\mathcal{D}_{\pi/2})
	= \mathbb{C} \backslash \{ y\mathrm{i} \mid |y| \geq \frac{\pi}{2} \}.
\end{align}
DE transformations have advantages that
they can avoid the singularities in $S_\psi$
and the parameter $C$ is large.
However, they have a significant problem that
they cannot avoid the singularities in $\tilde{S}$.
Specifically, the parameter $d$ may be quite small
if the integrand $f$ has singularities near the real axis.

%%%%%%%%%%
\subsection{Method of Slevinsky and Olver}
\label{subsec:SO}
%%%%%%%%%%

Slevinsky and Olver \cite{SO}
improved the DE formulas
by making polynomial adjustments to them.
They used the following formula as an option of the function $H$:
\begin{align}
	H_{\mathrm{SO}}(t) = C \sinh(t)
	+ \sum_{k=1}^m u_k t^{k-1},
\end{align}
where $C > 0$ and $u_1, \dots, u_m \in \mathbb{R}$.
Then, they chose these parameters
by solving an optimization problem.
In it, they maximized $C$
under the condition that the image of the boundary
$\partial \mathscr{D}_{\pi/2}$ under $H_{\mathrm{SO}}$
could pass through the singularities
$\{ \tilde{\delta}_j \pm \mathrm{i}\tilde{\epsilon}_j \}_{j=1}^m$.
This is formulated as follows:
\begin{align}
			\max &\quad C \quad 
			\text{ with respect to }
			C>0,\, u_1, \ldots, u_m,\, x_1, \ldots, x_m \in \mathbb{R} \\
			\label{eq:OptSO}
			\text{subject to} &\quad
			H_{\mathrm{SO}}\!\left(x_j + \frac{\pi}{2}\mathrm{i}\right)
			= \tilde{\delta}_j + \tilde{\epsilon}_j \mathrm{i}
			\quad (j=1, \ldots, m).
\end{align}
Algorithms to solve this is shown in \cite{SOprogram}.

However, we find some points to be improved
in their methods as follows:
\begin{itemize}
	\item
	They does not consider the singularities $S_\psi$.
	Thus, the parameters $C$ and $d$ may be smaller
	than the optimal values.
	\item
	The limiting conditions of the optimization problem \cref{eq:OptSO}
	does not imply the condition of the singularities \cref{eq:sing}.
	There are cases where $H_{\mathrm{SO}}$ is not an injection
	and $d < \pi/2$ even though
	the conditions \cref{eq:OptSO} are satisfied.
	\item
	There are cases where we cannot solve the optimization problem
	by \cite{SOprogram} because of the difficulty of
	finding the solution numerically.
	In these cases, we cannot use their methods.
\end{itemize}
We show experiments which cause the second and third problem
in \Cref{subsec:53,subsec:54}, respectively. 

%%%%%%%%%%%%%%%%%%%%%%%%%%%%%%%%%%%%%%%%%%%%%%%%%%%%%%%%%%%%%%%%%%%%%%%%%%%
\section{Proposed method}
\label{sec:4}
%%%%%%%%%%%%%%%%%%%%%%%%%%%%%%%%%%%%%%%%%%%%%%%%%%%%%%%%%%%%%%%%%%%%%%%%%%%

In this section, we propose a new method to construct the function $H$.
In it, we use a conformal map
from the domain $\mathscr{D}_{\pi/2}$ to a polygon $P$ 
as an option of $H$.
Then, choosing the function $H$ corresponds to choosing the polygon $P$.
It enables us to handle the condition of the singularities \cref{eq:sing} directly.
Also, using this conformal map has an advantage that it is written explicitly
using the idea of the Schwarz-Christoffel transformation.

%%%%%%%%%%
\subsection{Schwarz-Christoffel transformation}
\label{subsec:SC}
%%%%%%%%%%

\begin{figure}[ht]
\centering
\includegraphics[scale = 0.25]{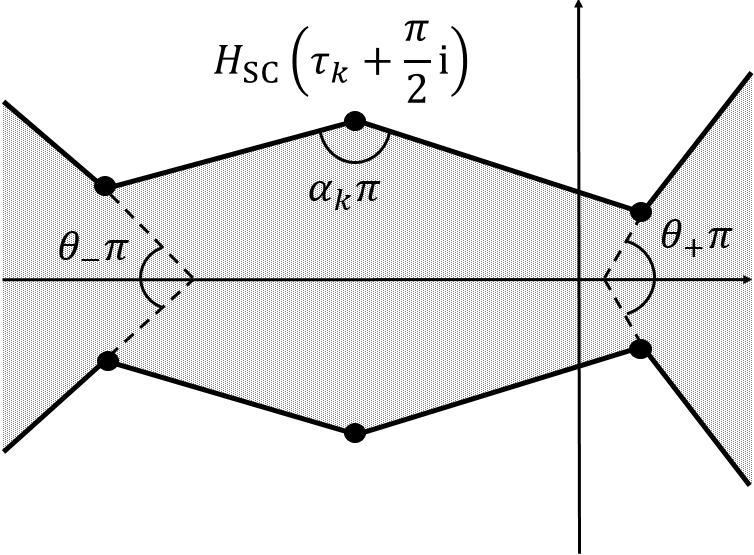}
\caption{Example of the polygon $P$. \label{img:SC}}
\end{figure}

The conformal mapping $H$ is written explicitly using
a modified version of the Schwarz-Christoffel transformation \cite{modSC}.
We assume that the polygon $P$ is symmetric with respect to the real axis
and that it has vertices at $\pm \infty$.
Let $P$ have $M$ pairs of the vertices other than $\pm \infty$ and
let $\{ \alpha_j \pi \}_{j=1}^M$ denote the interior angles of them.
We assume that $0 < \alpha_j \leq 2$ if the corresponding vertex is finite
and that $-2 \leq \alpha_j < 0$ if it is infinite.
Let $\theta_+$ and $\theta_-$ denote
the divergence angles of $\pm \infty$.
We assume that $0 \leq \theta_+, \theta_- \leq 1$.
\Cref{img:SC} shows an example of the polygon $P$.
Then, the conformal mapping from
the domain $\mathscr{D}_{\pi/2}$ to the polygon $P$
is written as the following theorem.

\begin{theorem} \label{thm:SC}
	For a given polygon $P$,
	there are some real numbers $\tau_1 < \tau_2 \cdots < \tau_M$
	such that a mapping
	\begin{align} \label{eq:HSC}
		H_{\mathrm{SC}}(z) =
		C\int_0^z
		\exp\left( \frac{1}{2} (\theta_{+}-\theta_{-}) \zeta \right)
		\left\{
			\prod_{j=1}^M
			\cosh^{\alpha_j-1} ( \zeta-\tau_j )
			\right\}
		\mathrm{d}\zeta + D
	\end{align}
	is a conformal mapping
	from the domain $\mathscr{D}_{\pi/2}$ to the polygon $P$
	which satisfies $H_{\mathrm{SC}}(+\infty) = +\infty$
	and $H_{\mathrm{SC}}(-\infty) = -\infty$.
	\begin{proof}
		A conformal mapping from the strip region
		$\{\zeta \mid 0 < \mathrm{Im}[\zeta] < 1\}$ to the polygon $P$
		has been shown in \cite{modSC}.
		We obtain \cref{eq:HSC}
		by transforming the domain linearly.
	\end{proof}
\end{theorem}
The problem of how to choose $\tau_1, \ldots, \tau_m$ is known
as the Schwarz-Christoffel parameter problem,
which has been studied in \cite{modSC,Trefethen}.

From \Cref{thm:SC}, 
the problem of how to choose the function $H$ is changed into the problem
of how to choose the polygon $P$.
We discuss it in the following subsection.

%%%%%%%%%%
\subsection{Suitable polygon for numerical integration}
\label{subsec:proposed}
%%%%%%%%%%

We consider how to construct the suitable polygon for numerical integration
according to the discussions in \Cref{subsec:H}.
The condition of the singularities \cref{eq:sing} is written simply as
\begin{align} 
	P \cap ( \tilde{S} \cup S_\psi ) = \emptyset .
\end{align}
We choose the polygon $P$ so that the corresponding parameters
$\gamma$ and $C$ will be larger under this condition.
For the following discussions, we rewrite
\begin{align}
	\tilde{S} \cup S_\psi = 
	\{ \tilde{\delta}_j \!\pm\! (\tilde{\epsilon}_j +2k\pi)\mathrm{i}
	\mid j = 1, \ldots, m, k\in \mathbb{Z} \},
\end{align}
where $\tilde{\delta}_1 < \dots < \tilde{\delta}_m$ and
$0 < \tilde{\epsilon}_j \leq \pi\,(j = 1, \ldots, m)$.

Relations between the polygon $P$ and the parameter $\gamma$ are obtained
by asymptotic expansion of $H_{\mathrm{SC}}$.
\begin{theorem} \label{thm:asympt}
	Let $\{\alpha_j\}_{j=1}^M, \theta_+$, and $\theta_-$ be the parameters
	introduced in \Cref{subsec:SC}. We write
	$\bar{\theta} = (\theta_+ + \theta_-)/2$
	and $\Delta\theta = (\theta_+ - \theta_-)/2$.
	Then,
	\begin{align} \label{eq:asympt}
		& \int_0^t \mathrm{e}^{\Delta\theta \tau}
		\prod_{j=1}^M \cosh^{\alpha_j -1}(\tau-\tau_j)
		\mathrm{d}\tau \\
		=&\ 
		\frac{1}{\theta_+ \theta_-2^{\bar\theta-1}}
		\frac{1}{2}
		\left(
			\theta_-
			\mathrm{e}^{\theta_+ t- \sum_{j=1}^M (\alpha_j - 1) \tau_j}
			- \theta_+
			\mathrm{e}^{-\theta_- t+ \sum_{j=1}^M (\alpha_j - 1) \tau_j}
		\right)
		+ \mathrm{O}(1) \nonumber
	\end{align}
	holds as $|t| \to \infty, t \in \mathbb{R}$.
	Specifically, when the divergence angles of the both sides are equal,
	i.e., $\theta_+ = \theta_- = \bar{\theta}$,
	\begin{align} \label{eq:asympt_symmetry}
		\int_0^t
		\prod_{j=1}^M \cosh^{\alpha_j -1}(\tau-\tau_j)
		\mathrm{d}\tau 
		= \frac{1}{2^{\bar{\theta}-1}\bar{\theta}^2}
		\sinh(\bar{\theta}t - \sum_{j-1}^M (\alpha_j - 1) \tau_j)
		+ \mathrm{O}(1)
	\end{align}
	holds as $|t| \to \infty, t \in \mathbb{R}$.
\end{theorem}
The proof is given in \Cref{ap:asympt}.

From this theorem,
we see that $\gamma = \min \{\theta_+, \theta_- \}$.
Thus, we fix $\theta_+$ and $\theta_- $ to 1.

Relations between the polygon $P$ and the parameter $C$
are rather complex.
Although it is difficult to formulate them,
it is observed experimentally that
the larger the area of the polygon $P$ is,
the larger the parameter $C$ is.
We show the experiments in \Cref{ap:C}.

\begin{figure}[ht]
	\centering
	\includegraphics[scale = 0.25]{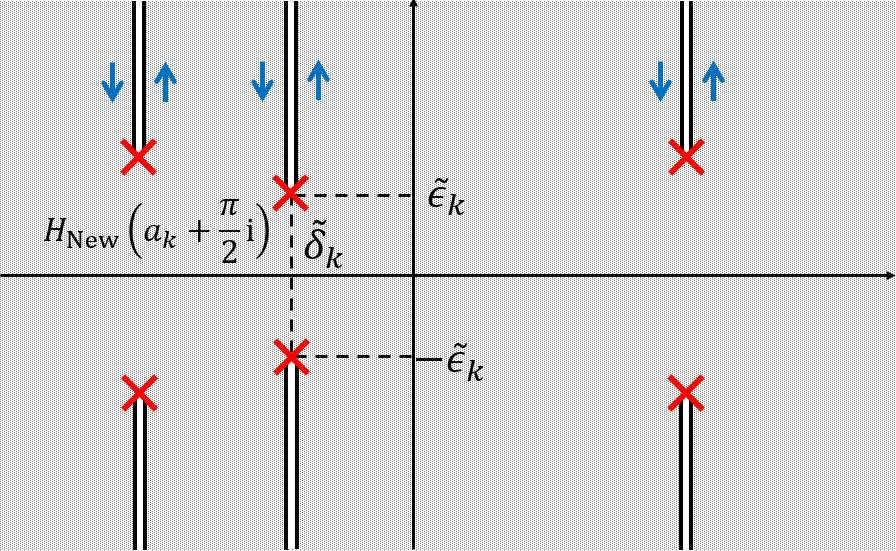}
	\caption{The proposed polygon $P_{\mathrm{New}}$.
		\label{img:proposed}}
\end{figure}

For these reasons,
we propose to construct a new transformation formula $H_{\mathrm{New}}$
so that the image $H_{\mathrm{New}} (\mathscr{D}_{\pi/2})$ will avoid
singularities by $m$ pairs of slits as \Cref{img:proposed}.
Let $P_{\mathrm{New}}$ denote this.
The function $H_{\mathrm{New}}$ coincides with
DE transformations if $\tilde{S} = \emptyset$ and
$S_\psi = \{ ( \pm \pi/2 + 2k\pi)\mathrm{i} \mid k \in \mathbb{Z} \}$.
It can be considered
to be a generalization of DE transformations.

The corresponding transformation to $P_{\mathrm{New}}$
is written as follows.
The polygon $P_{\mathrm{New}}$ has $(2m-1)$ pairs of vertices
other than $\pm \infty$
at the turning points of the slits (with angles $2\pi$)
and the points at infinity (with angles 0).
Let $\{ a_j \pm \frac{1}{2}\pi\mathrm{i} \}_{j=1}^m$
and $\{ b_j \pm \frac{1}{2}\pi\mathrm{i} \}_{j=1}^{m-1}$ denote
the preimages of these vertices under $H_{\mathrm{New}}$, respectively.
Then, we can write $(\tau_1, \ldots, \tau_{2m-1}) = (a_1, b_1, a_2, \ldots, b_{m-1}, a_m)$
and $(\alpha_1, \ldots, \alpha_{2m-1}) = (2,0,2,\ldots,0,2)$.
The transformation $H_{\mathrm{New}}$ is obtained by substituting these parameters
into \cref{eq:HSC}, that is,
\begin{align} \label{eq:HNew_0}
	H_{\mathrm{New}}(z)
	= C\int_0^z
	\frac
		{\prod_{j=1}^m\cosh(z-a_j)}
		{\prod_{j=1}^{m-1}\cosh(z-b_j)}
	\mathrm{d}z + D,
\end{align}
where $C>0$ and $a_1<b_1<\cdots<b_{m-1}<a_m$.
Using \Cref{thm:asympt}, the asymptotic form of $H_{\mathrm{New}}$ is written as
\begin{align} \label{eq:HNew_asympt}
	H_{\mathrm{New}}(t) 
	&= C\sinh(t-T) + \mathrm{O}(1)
	\quad ( |t| \to \infty,\, t \in \mathbb{R}),
\end{align}
where $T = a_1 - b_1 + a_2 - \cdots - b_{m-1} + a_m$.

The parameters are determined as follows.
First, we can choose the parameter $T$ arbitrarily.
We determine it so that the parameter $\beta_2$ will be the largest.
Then, we determine the other parameters
$C, a_1, \ldots, a_m$ and $b_1, \ldots, b_{m-1}$
so that the image
$H_{\mathrm{New}}(\mathscr{D}_{\pi/2})$ will match $P_{\mathrm{New}}$.
It has known that these parameters
are uniquely determined with the value of $T$ fixed \cite{modSC}.
We show how to calculate these parameters
in \Cref{subsec:T,subsec:bCD}, respectively.

Incidentally, there is another advantage to construct the transformation
$H_{\mathrm{New}}$ in this way.
The formula $H_{\mathrm{New}}$ is written rather simply.
This is important when we use $H_{\mathrm{New}}$ as a change of variables.
We show this in \Cref{subsec:bCD}.

%%%%%%%%%%
\subsection{Determination of the parameter \textit{T}}
\label{subsec:T}
%%%%%%%%%%

Here, we consider an integral of the interval $(-1, 1)$.
The cases of the other intervals are shown
in \Cref{ap:T}.

We assume that the integrand $f$ is smooth on the interval $(-1, 1)$
and satisfies
\begin{align}
	\label{orderf_finite}
	f(x) = 
	\left\{
		\begin{aligned}
			&\mathrm{O}( (1-x)^p ) \quad ( x \to 1-0) \\
			&\mathrm{O}( (1+x)^q ) \quad ( x \to -1+0)
		\end{aligned}
	\right.
\end{align}
for some $p, q > -1$.
The change of variables is given as
\begin{align}
	x = \phi(t) = \tanh(H_{\mathrm{New}}(t)).
\end{align}
The decay rate of the transformed integrand is estimated as
\begin{align}
	f(\phi(t))\phi'(t) = 
	\left\{
		\begin{aligned}
			&\mathrm{O}\!\left(\exp( -(C(1+p) -\varepsilon) \mathrm{e}^{t-T}) \right)
			\quad (t\to+\infty) \\
			&\mathrm{O}\!\left(\exp( -(C(1+q)-\varepsilon) \mathrm{e}^{T-t}) \right)
			\quad (t\to-\infty)
		\end{aligned}
	\right.
\end{align}
for arbitrary $\varepsilon >0$.
Then we see that the parameter $\beta_2$ satisfies
\begin{align}
	\beta_2 \leq
	\min\!
	\left\{
		(C(1+p) -\varepsilon)\mathrm{e}^{-T},\, 
		(C(1+q) -\varepsilon)\mathrm{e}^{T}
	\right\} .
\end{align}

To make the parameter $\beta_2$ larger,
we we make $\varepsilon$ go to 0 and determine $T$ as
\begin{align}
	-\frac{C}{2}(1+r)\mathrm{e}^{-T} &= \frac{C}{2}(1+q)\mathrm{e}^T
	\quad \Leftrightarrow \quad 
	T = \frac{1}{2}
	\log\! \left(\! - \frac{1+r}{1+q} \right).
\end{align}
Then the supremum of the parameter $\beta_2$ is estimated as
\begin{align}
	\beta_2^* = C\sqrt{(p+1)(q+1)}.
\end{align}

%%%%%%%%%%
\subsection{Determination of the other parameters}
\label{subsec:bCD}
%%%%%%%%%%

First, we discuss the case of $m = 2$:
\begin{align} \label{eq:HNew_case2_0}
	H_{\mathrm{New}}(z)
	= C\int_0^z
	\frac{\cosh(z-a_1)\cosh(z-a_2)}{\cosh(z-b_1)}
	\mathrm{d}z + D \quad  (C>0,\, a_1<b_1<a_2).
\end{align}
The integrand of \cref{eq:HNew_case2_0} is rearranged as
\begin{align}
	\frac{\cosh(z-a_1)\cosh(z-a_2)}{\cosh(z-b_1)}
	=&\ 
	\frac{1}{2}
	\frac
	{\mathrm{e}^{-a_1+b_1-a_2}\left(\mathrm{e}^{2z}+\mathrm{e}^{2a_1}\right)
	\left(\mathrm{e}^{2z}+\mathrm{e}^{2a_2}\right)}
	{\mathrm{e}^z\left(\mathrm{e}^{2z}+\mathrm{e}^{2b_1}\right)} \\
	=&\ 
	\cosh(z-a_1+b_1-a_2)+
	\frac{L_1\mathrm{e}^{z-b_1}}{\mathrm{e}^{z-b_1}+\mathrm{e}^{-z+b_1}},
\end{align}
where
$ 2L_1=
\mathrm{e}^{-a_1+b_1-a_2}\!
\left(\mathrm{e}^{2a_1}+\mathrm{e}^{2a_2}-\mathrm{e}^{2b_1}
-\mathrm{e}^{2\left(a_1-b_1+a_2\right)}\right)$.
Then, the transformation $H_{\mathrm{New}}$ is written as
\begin{align}
	H_{\mathrm{New}}(z)
	&=
	C\sinh(z-a_1+b_1-a_2)+
	CL_1\int_0^z
	\frac
	{\mathrm{e}^{z-b_1}}
	{\mathrm{e}^{z-b_1}+\mathrm{e}^{-z+b_1}} + D \\
	&=
	C\sinh(z-T)
	+ 2D_1\tan^{-1}\!\left(\mathrm{e}^{z-b_1} \right)
	+ D_0,
\end{align}
where
$T = a_1-b_1+a_2$,
$2D_1 = CL_1$, and $D_0 = D-CL_1\tan^{-1}\!\left(\mathrm{e}^{-b_1}\right)$.
The value of the parameter $T$ has been determined in \Cref{subsec:T}.
We determine the other parameters
so that
the image of the upper boundary of the strip region $\mathscr{D}_{\pi/2}$
will match the upper slits of \Cref{img:proposed}.
For this reason, we consider the image of
$z = x + \frac{\pi}{2}\mathrm{i}, x \in \mathbb{R}$:
\begin{align} \label{eq:HNew_case2_up}
	H_{\mathrm{New}}\!\left(x+\frac{\pi}{2}\mathrm{i}\right)
	= C\cosh(x-T)\mathrm{i}
	-D_1 \log\!
	\left( \frac {1-\mathrm{e}^{x-b_1}}{1+\mathrm{e}^{x-b_1}} \right)
	\mathrm{i}+D_0.
\end{align}

The real part of \cref{eq:HNew_case2_up} is given by
\begin{align}
	\mathrm{Re} \left[H_{\mathrm{New}}\!\left(x+\frac{\pi}{2}\mathrm{i}\right)\right]
	=
	D_1 \mathrm{arg} \!
	\left(\frac{1-\mathrm{e}^{x-b_1}}{1+\mathrm{e}^{x-b_1}}
	\right)
	+D_0
	= \left\{
	\begin{aligned}
		&D_0 & & (x < b_1) \\
		&D_0 + \pi D_1 & & (x > b_1).
	\end{aligned} \right.
\end{align}
Thus, we determine $D_0$ and $D_1$ as
\begin{align}
	\left\{
	\begin{aligned}
		D_0 &= \tilde{\delta}_1 \\
		D_0 + \pi D_1 &= \tilde{\delta}_2
	\end{aligned} \right.
	\quad\Leftrightarrow\quad
	\left\{
	\begin{aligned}
		D_0 &= \tilde{\delta}_1 \\
		D_1 &= \frac{1}{\pi}\!\left(\tilde{\delta}_2-\tilde{\delta}_1\right).
	\end{aligned} \right.
\end{align}
The imaginary part of \cref{eq:HNew_case2_up} is given by
\begin{align} \label{eq:HNew_case2_up_imag}
	\mathrm{Im}
	\left[ H_{\mathrm{New}}\!\left(x+\frac{\pi}{2}\mathrm{i}\right) \right] 
	= C\cosh(x-T) - D_1\log
	\left| \tanh \left(\frac{x-b_1}{2}\right) \right|.
\end{align}
The function \cref{eq:HNew_case2_up_imag} has local minimum points
in $( -\infty, b_1 )$ and $( b_1, \infty )$, which correspond to
the parameters $a_1$ and $a_2$, respectively.
The function values at them correspond to
$\tilde{\epsilon}_1$ and $\tilde{\epsilon}_2$.
Thus, we determine the parameters $C, a_1, b_1$, and $a_2$ by solving 
\begin{align}
	\label{eq:HNew_case2_solve}
	\left\{
		\begin{aligned}
			C\cosh(a_1-T)
			- D_1\log\!\left|\tanh (a_1-b_1) / 2 \right|
			&= \tilde{\epsilon}_1 \\
			C\cosh(a_2-T)
			- D_1\log\!\left|\tanh (a_2-b_1) / 2 \right|
			&= \tilde{\epsilon}_2 \\
			C\sinh(a_1-T)
			- D_1 / \sinh(a_1-b_1)
			&= 0 \\
			C\sinh(a_2-T)
			- D_1 / \sinh(a_2-b_1)
			&= 0 \\
		\end{aligned}
	\right.
\end{align}
under the constraints $C > 0$ and $a_1 < b_1 < a_2$.
Here, the condition $T = a_1 - b_1 + a_2$ is automatically satisfied
by solving \cref{eq:HNew_case2_solve}.
We show it later in the general case.

Then, we extend the discussions to the general case.
The following proposition shows
that we can deform the integrand similarly to the case of $m=2$.
\begin{proposition}
	We write $T = a_1-b_1+\cdots-b_{m-1}+a_m$. Then,
	\begin{align} \label{eq:Hexpand}
		\frac{\prod_{j=1}^m\cosh(z-a_j)}{\prod_{j=1}^{m-1}\cosh(z-b_j)}
		&\equiv
		\cosh(z-T)+
		\sum_{j=1}^{m-1}
		\frac
		{L_j\mathrm{e}^{z-b_j}}
		{\mathrm{e}^{z-b_j}+\mathrm{e}^{-z+b_j}}
	\end{align}
holds for some $L_1, \ldots, L_{m-1} \in \mathbb{R}$.
	\begin{proof}
		We define $m$th degree polynomials $F_1$ and $F_2$ as
		$F_1(Z) = \prod_{j=1}^m (Z+\mathrm{e}^{2a_j})$ and
		$F_2(Z) = (Z+\mathrm{e}^{2T})\prod_{j=1}^{m-1} (Z+\mathrm{e}^{2b_j})$.
		Since the leading terms and constants of $F_1$ and $F_2$ coincide,
		We can write
		\begin{align}
			F_1(Z) = F_2(Z) + ZG(Z)
		\end{align}
		for some polynomial $G$ of which the degree is $(m-2)$ or less.
		Also, using Lagrange polynomials, the polynomial $G$ is written as
		\begin{align}
			G(Z) = \sum_{j=1}^{m-1} G(-\mathrm{e}^{2b_j})
		\frac
		{\prod_{k=1, \dots, m-1, k\neq j} (Z+\mathrm{e}^{2b_k})}
		{\prod_{k=1, \dots, m-1, k\neq j} (-\mathrm{e}^{2b_j}+\mathrm{e}^{2b_k})}
		=:\sum_{j=1}^{m-1}l_j\prod_{k\neq j} (Z+\mathrm{e}^{2b_k}).
		\end{align}
		Then, we obtain \cref{eq:Hexpand} by
		\begin{align}
			\frac{\prod_{j=1}^m\cosh(z-a_j)}{\prod_{j=1}^{m-1}\cosh(z-b_j)}
			&= \frac{1}{2}
			\frac{\mathrm{e}^{-T}F_1 (\mathrm{e}^{2z})}
			{\mathrm{e}^{z}\prod_{j=1}^{m-1} (\mathrm{e}^{2z}+\mathrm{e}^{2b_j})} \\
			&= \frac{1}{2}
			\frac
			{\mathrm{e}^{-T}\left\{F_2 (\mathrm{e}^{2z})
				+ \mathrm{e}^{2z}G (\mathrm{e}^{2z}) \right\}}
			{\mathrm{e}^{z}\prod_{j=1}^{m-1}
				(\mathrm{e}^{2z}+\mathrm{e}^{2b_j})} \\
			&=
			\cosh(z-T) + \sum_{j=1}^{m-1}
			\frac
			{L_j\mathrm{e}^{z-b_j}}
			{\mathrm{e}^{z-b_j}+\mathrm{e}^{-z+b_j}},
		\end{align}
		where $L_j = \mathrm{e}^{-T}l_j / 2$.
	\end{proof}
\end{proposition}
Using this proposition, we can rewrite $H_{\mathrm{New}}$ similarly as
\begin{align}
	H_{\mathrm{New}}(z) = C\sinh(z-T)
	+ \sum_{j=1}^{m-1} 2D_j\tan^{-1}\!\left(\mathrm{e}^{z-b_j}\right)+D_0
\end{align}
for some $D_0, D_1, \ldots, D_m \in \mathbb{R}$.

The parameters are also determined similarly to the case of $m = 2$.
The parameter $T$ has been determined in \Cref{subsec:T}.
The parameters $D_0, D_1, \ldots, D_{m-1}$ are determined
based on the real part of $H_{\mathrm{New}}\!\left(x+\frac{\pi}{2}\mathrm{i}\right)$: 
\begin{align}
	\left\{
		\begin{aligned}
			D_0 &=\tilde{\delta}_1 \\
			D_j &=\frac{1}{\pi}\!\left(\tilde{\delta}_{j+1}-\tilde{\delta}_j\right)
			\quad\left(j=1,\, \ldots,\, m-1\right).
		\end{aligned}
	\right.
\end{align}
The other parameters are determined
based on the imaginary part of $H_{\mathrm{New}}(x+(\pi/2)\mathrm{i})$,
that is,
the parameters $C$, $a_1, \ldots, a_m$ and $b_1, \ldots, b_{m-1}$
are determined by solving
\begin{align} \label{eq:HNew_solve}
	\left\{
		\begin{aligned}
			C\cosh(a_k-T)
			-\sum_{j=1}^{m-1}
			D_j\log\!\left|\tanh\!\left(\frac{a_k-b_j}{2}\right)\right|
			&=\tilde{\epsilon}_k
			&\quad&
			(k=1,\, \ldots,\, m) \\
			C\sinh(a_k-T)
			-\sum_{j=1}^{m-1}\frac{D_j}{\sinh(a_k-b_j)}
			&=0
			&\quad&
			(k=1,\, \ldots,\, m)
		\end{aligned}
	\right.
\end{align}
under the constraints $C>0$ and $a_1<b_1<\cdots<b_{m-1}<a_m$.
Here, the condition $T = a_1 - b_1 + \cdots - b_{m-1} + a_m$
is automatically satisfied by solving \cref{eq:HNew_solve}.
The following theorem shows this.
\begin{theorem} \label{thm:ababa}
	We assume that a system of equations
	\begin{align}
		\left\{
			\begin{aligned}
				C\sinh(a_1-T)\ - &\sum_{j=1}^{m-1}
				\frac{D_j}{\sinh(a_1-b_j)} = 0 \\
				&\quad\vdots \\
				C\sinh(a_m-T)\ - &\sum_{j=1}^{m-1}
				\frac{D_j}{\sinh(a_m-b_j)} = 0
			\end{aligned}
		\right.
	\end{align}
	holds for some real numbers $a_1, \ldots, a_m$,
	$b_1, \ldots, b_{m-1}$, $D_0, \ldots, D_{m-1}$, C and T
	which satisfy
	$a_1<b_1<\cdots<b_{m-1}<a_m$ and $(C, D_1, \ldots, D_m) \neq (0, \ldots, 0)$.
	Then,
	\begin{align}
		a_1 - b_1 + \cdots - b_{m-1} + a_m = T
	\end{align}
	holds.
\end{theorem}
The proof is given in \Cref{ap:ababa}.

In general, it is difficult to solve a non-linear system of equations such as
\cref{eq:HNew_solve} numerically under constraints.
Thus, we replace the parameters with $x \in \mathbb{R}^{2m}$ as
\begin{align}
	\begin{aligned}
		\left\{
		\begin{aligned}
			x_1 &=\log C \\
			x_2 &=a_1 \\ 
			x_3 &=\log\!\left(b_1-a_1\right) \\
			x_4 &=\log\!\left(a_2-b_1\right) \\
			x_5 &=\log\!\left(b_2-a_2\right) \\
			\vdots \\
			x_{2m}&=\log\!\left(a_{m}-b_{m-1}\right)
		\end{aligned}
		\right.
		\quad \Leftrightarrow \quad
		\begin{aligned}
			\left\{
			\begin{aligned}
				C &=\mathrm{e}^{x_1} \\
				a_1 &= x_2 \\
				b_1 &= x_2+\mathrm{e}^{x_3} \\
				a_2 &= x_2+\mathrm{e}^{x_3}+\mathrm{e}^{x_4}\\
				\vdots \\
				a_m &= x_2+\mathrm{e}^{x_3}+\mathrm{e}^{x_4}
				+\ldots+\mathrm{e}^{x_{2m}},
				\end{aligned}
			\right.
		\end{aligned}
	\end{aligned}
\end{align}
and we solve \cref{eq:HNew_solve} as a system of equations of $x$.
This is a similar method to numerical computations
of Schwarz-Christoffel transformation \cite{Trefethen}.
We solve it using NLsolve,
a Julia program for solving non-linear systems of equations \cite{NLsolve}.

%%%%%%%%%%
\subsection{Approximation of the proposed transformation}
\label{subsec:4_approx}
%%%%%%%%%%

The proposed transformation $H_{\mathrm{New}}$ has a disadvantage that
calculating the terms of $\tan^{-1}$ takes a lot of time.
For this reason, we also propose approximating them by
\begin{align} \label{eq:approx}
	\tan^{-1}\!\left(\mathrm{e}^t\right)
	\approx \frac{\pi}{4}\!
	\left(\tanh\!\left(\frac{2}{\pi}t \right) +1 \right),
\end{align}
where the values and the derivatives at $t = 0$,
and the limits as $t \to \pm \infty$ of both sides coincide.
We construct an approximation formula $H_{\mathrm{New2}}$
by approximating $H_{\mathrm{New}}$ using \cref{eq:approx}:
\begin{align}
		H_{\mathrm{New2}}(t) = C\sinh(t-T)
		+ \left\{ \sum_{j=1}^{m-1}
		\frac{\pi}{2} D_j
		\left( \tanh\!\left( \frac{2}{\pi}(t-b_j)\! \right)\! +1 \right)
		\right\} + \tilde{\delta}_1,
\end{align}
where the parameters are determined
by the methods of Subsection \Cref{subsec:T,subsec:bCD}.
%
%%%%%%%%%%%%%%%%%%%%%%%%%%%%%%%%%%%%%%%%%%%%%%%%%%%%%%%%%%%%%%%%%%%%%%%%%%%
\section{Numerical experiments}
\label{sec:5}
%%%%%%%%%%%%%%%%%%%%%%%%%%%%%%%%%%%%%%%%%%%%%%%%%%%%%%%%%%%%%%%%%%%%%%%%%%%

We compare the effectiveness of the transformation formulas $H$
by some numerical experiments.
In \Cref{subsec:51,subsec:52},
we deal with the same examples as those in \cite{SO}.
In \Cref{subsec:53},
we show an example
that the transformation $H_{\mathrm{SO}}$ is not an injection.
In \Cref{subsec:54},
we show an example 
to which we cannot apply the method of Slevinsky and Olver.

%%%%%%%%%%
\subsection{Integral on a finite interval}
\label{subsec:51}
%%%%%%%%%%

We consider an integral on a finite interval \cite{SO}:
\begin{align}
	\int_{-1}^{1}
	\frac
	{\exp \left( (\epsilon_1^2 + (x-\delta_1)^2 )^{-1}\right) \log(1-x)}
	{(\epsilon_2^2 + (x-\delta_2)^2 ) \sqrt{1+x}}
	\mathrm{d}x = -2.04645 \ldots ,
\end{align}
where $\delta_1 \pm \epsilon_1\mathrm{i} = - 0.5 \pm \mathrm{i}$ and
$\delta_2 \pm \epsilon_2\mathrm{i} = - 0.5 \pm 0.5\mathrm{i}$.
The change of variables for this integral is $\phi (t) = \tanh (H(t))$.
There are singularities as follows:
\begin{align}
	\tilde{S} = \{ \tanh^{-1}(\delta_j \pm \epsilon_j \mathrm{i} )
	\mid j = 1, 2,\, k \in \mathbb{Z} \}, \quad
	S_\psi = \left\{ 
	\left(\pm \frac{1}{2} + 2k \right)\pi \mathrm{i} \mid k \in \mathbb{Z} \right\}.
\end{align}

\begin{figure}[ht]
	\centering
	\subfigure[$H_{\mathrm{SO}}(\mathscr{D}_{\frac{\pi}{2}})$]
	{\includegraphics[scale = 0.26]{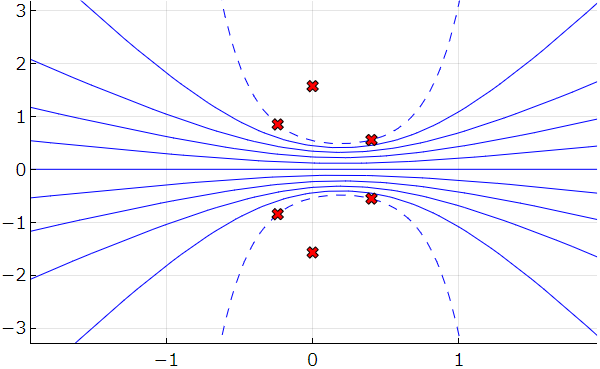}}
	\subfigure[$H_{\mathrm{New}}(\mathscr{D}_{\frac{\pi}{2}})$]
	{\includegraphics[scale = 0.26]{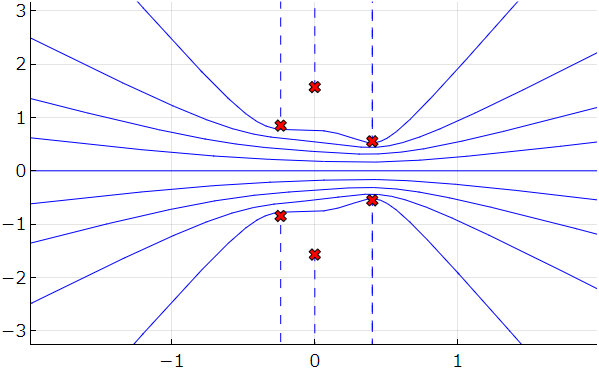}}
	\subfigure[$H_{\mathrm{New2}}(\mathscr{D}_{\frac{\pi}{2}})$]
	{\includegraphics[scale = 0.26]{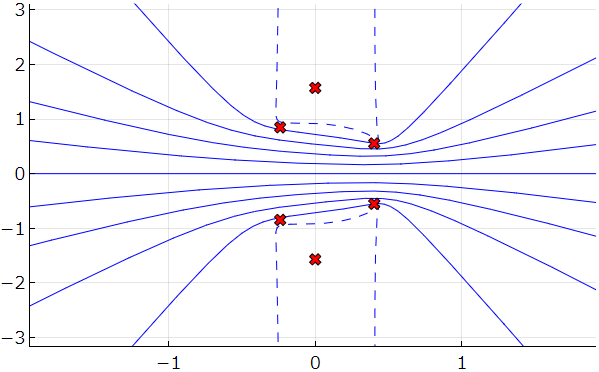}}
	\caption
	{Images $H(\mathscr{D}_{\pi/2})$ and singularities.
	The solid lines are the images of lines
	which are parallel to the real axis in $\mathscr{D}_{\pi/2}$.
	The dotted lines show $H(\partial \mathscr{D}_{\pi/2})$.
	\label{img:51H}}
\end{figure}

First, the formulas $H$ are given by
\begin{align}
	H_{\mathrm{DE}}(t) =&\ \frac{\pi}{2}\sinh(t), \\
	H_{\mathrm{SO}}(t) \approx&\ 0.139 \sinh(t) + 0.191 + 0.219\,t, \\
	H_{\mathrm{New}}(t) \approx&\ 0.356\sinh(t-0.347) \\
	&+ 0.152 \tan^{-1}\left(\mathrm{e}^{t + 0.190} \right)
	+ 0.256 \tan^{-1}\!\left(\mathrm{e}^{t + 0.177}\right) - 0.239, \nonumber\\
	H_{\mathrm{New2}}(t) \approx&\ 0.356\sinh(t-0.347) \\
	&+ 0.119 \tanh(t + 0.190) + 0.201 \tanh(t + 0.177) - 0.0817 \nonumber.
\end{align}
\Cref{img:51H} shows the images $H(\mathscr{D}_{\pi/2})$.

\begin{figure}
	\centering
	\subfigure[]
	{\includegraphics[scale = 0.4]{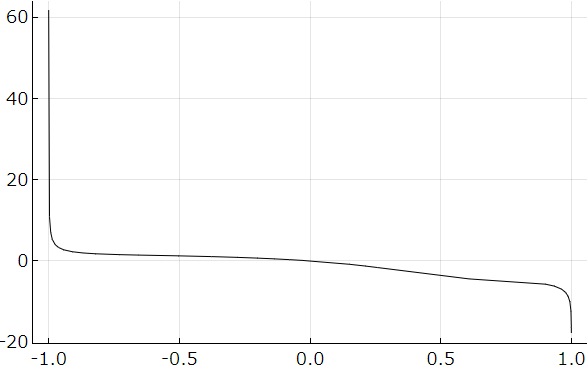}}
	\subfigure[]
	{\includegraphics[scale = 0.4]{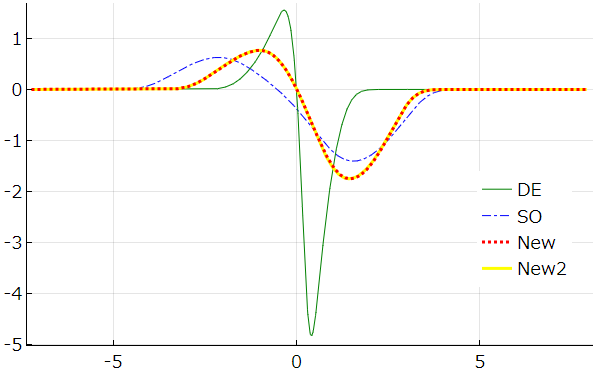}}
	\caption
	{{\rm (a)} The original integrand $f$.
	{\rm (b)} The transformed integrands $f(\phi(\cdot))\phi '(\cdot)$.
	\label{img:51F}}
\end{figure}

Then, we compare the performances of the formulas $H$
as transformation formulas for integration.
\Cref{table:51} shows the parameters of \Cref{thm:Sugihara_DE},
where the parameter $d$ of DE is calculated as 
\begin{align}
	d_{\mathrm{DE}}
	= \min_{j=1, 2}\ \mathrm{Im} \left[
	\sinh^{-1}\!\left(\frac{2}{\pi}\tanh^{-1}
	(\delta_j + \epsilon_j\mathrm{i})\right) \right] .
\end{align}
\Cref{img:51F} shows the original and transformed integrands.
\begin{table}[ht]
	\footnotesize
	\centering
	\caption{Parameters of \Cref{thm:Sugihara_DE}. \label{table:51}}
	\begin{tabular}{|c|cccc|}
		\hline
		& DE & SO & New & New2 \\
		\hline
		$\gamma$ & 1 & 1 & 1 & 1 \\
		$d$ & 0.346 & $\pi/2$ & $\pi/2$ & ?\\
		$\beta_2$ & 0.785 & 0.0695 & 0.252  & 0.252 \\
		\hline
	\end{tabular}
\end{table}

\begin{figure}
	\subfigure[\label{img:51err}]
	{\includegraphics[scale = 0.4]{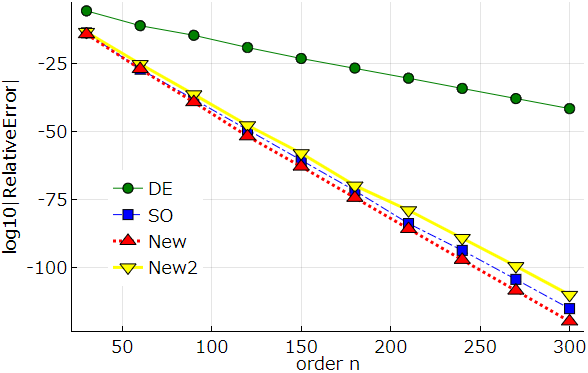}}
	\subfigure[\label{img:51time}]
	{\includegraphics[scale = 0.4]{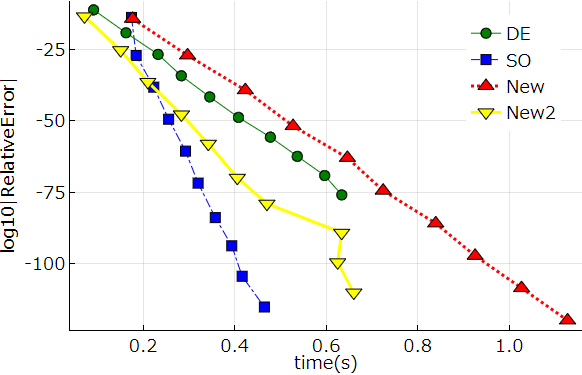}}
	\caption{{\rm (a)} Orders $n$ and errors.
	{\rm (b)} Calculation time and errors.
	The calculation time includes time
	to determine parameters and carry out the trapezoidal formula.}
\end{figure}

Finally, we compare the errors of the numerical integration.
\Cref{img:51err} shows relations between orders $n$ and
the errors. \Cref{img:51time} shows relations between time
for the calculation of numerical integration and the errors.
Here, we assume that $d_{\mathrm{New2}} = \pi/2$
when we calculate the mesh size of the trapezoidal formula \cref{eq:hDE}.

%%%%%%%%%%
\subsection{Integral on the infinite interval}
\label{subsec:52}
%%%%%%%%%%

We consider an integral on the infinite interval \cite{SO}:
\begin{align}
	\int_{-\infty}^{\infty}
	\frac
	{\exp\!\left(10(\epsilon_1^2 + (x-\delta_1)^2)^{-1} \right)\!
	\cos\!\left(10(\epsilon_2^2 + (x-\delta_2)^2)^{-1} \right)}
	{( (x-\delta_3)^2+\epsilon_3^2 )
	\sqrt{(x-\delta_4)^2+\epsilon_4^2}}
	\mathrm{d}x = 15.0136\ldots,
\end{align}
where
$\delta_1 \pm \epsilon_1\mathrm{i} = - 2 \pm \mathrm{i}$, 
$\delta_2 \pm \epsilon_2\mathrm{i} = - 1 \pm 0.5\mathrm{i}$, 
$\delta_3 \pm \epsilon_3\mathrm{i} = 1 \pm 0.25\mathrm{i}$, 
and $\delta_4 \pm \epsilon_4\mathrm{i} = 2 \pm \mathrm{i}$.
The change of variables for this integral is $\phi (t) = \sinh (H(t))$.
There are singularities as follows:
\begin{align}
	\tilde{S} = \{ \tanh^{-1}(\delta_j \pm \epsilon_j \mathrm{i} )
	\mid j = 1, \ldots, 4, \, k \in \mathbb{Z} \}, \quad
	S_\psi = \emptyset .
\end{align}

\begin{figure}[ht]
	\centering
	\subfigure[$H_{\mathrm{SO}}(\mathscr{D}_{\frac{\pi}{2}})$]
	{\includegraphics[scale = 0.26]{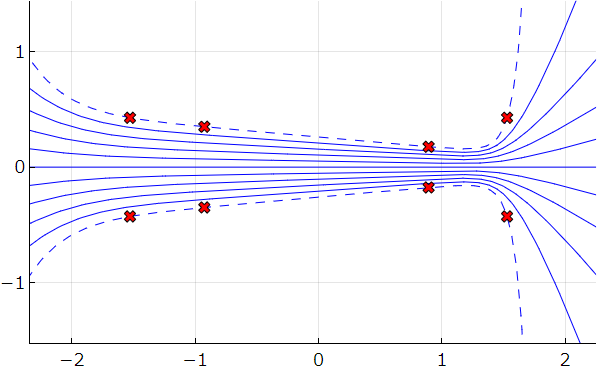}}
	\subfigure[$H_{\mathrm{New}}(\mathscr{D}_{\frac{\pi}{2}})$]
	{\includegraphics[scale = 0.26]{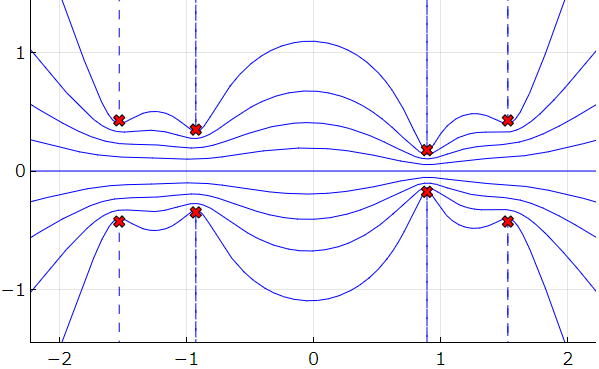}}
	\subfigure[$H_{\mathrm{New2}}(\mathscr{D}_{\frac{\pi}{2}})$]
	{\includegraphics[scale = 0.26]{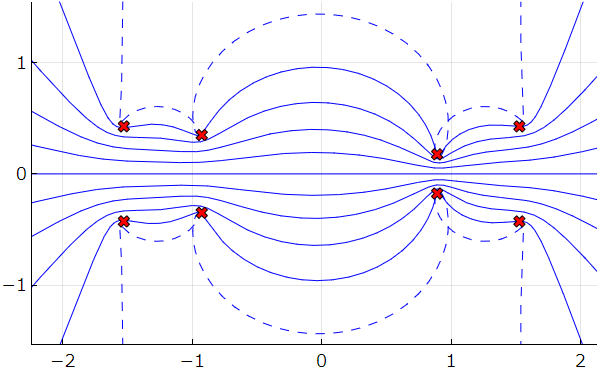}}
	\caption
	{Images $H(\mathscr{D}_{\pi/2})$ and singularities.
	The solid lines are the images of lines
	which are parallel to the real axis in $\mathscr{D}_{\pi/2}$.
	The dotted lines show $H(\partial \mathscr{D}_{\pi/2})$.
	\label{img:52H}}
\end{figure}

First, the formulas $H$ are given by
\begin{align}
	H_{\mathrm{DE}}(t) =&\ \frac{\pi}{2}\sinh(t), \\
	H_{\mathrm{SO}}(t) \approx&\ 
	5.77\times 10^{-6} \sinh(t) - 0.254 \\
	& + 0.149\,t - 4.54\times 10^{-3}\,t^2 + 9.99\times 10^{-5}\,t^3 , \nonumber \\
	H_{\mathrm{New}}(t) \approx&\ 5.12\times 10^{-3}\sinh(t)
	+ 0.384 \tan^{-1}\!\left(\mathrm{e}^{t + 4.32}\right) \\
	&+ 1.15 \tan^{-1}\!\left(\mathrm{e}^{t + 1.37}\right)
	+ 0.405 \tan^{-1}\!\left(\mathrm{e}^{t - 2.98}\right) - 1.53, \nonumber \\
	H_{\mathrm{New2}}(t) \approx&\ 5.12\times 10^{-3}\sinh(t) 
	+ 0.309 \tanh(t + 4.32) \\
	& + 0.909 \tanh(t + 1.37) + 0.318 \tanh(t - 2.98) \nonumber.
\end{align}
\Cref{img:52H} shows the images $H(\mathscr{D}_{\pi/2})$.

\begin{figure}
	\centering
	\subfigure[]
	{\includegraphics[scale = 0.39]{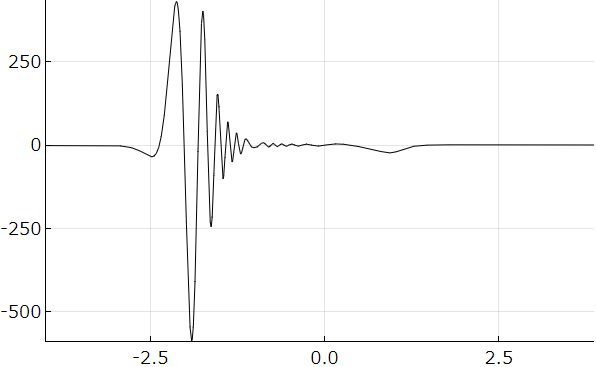}}
	\subfigure[]
	{\includegraphics[scale = 0.39]{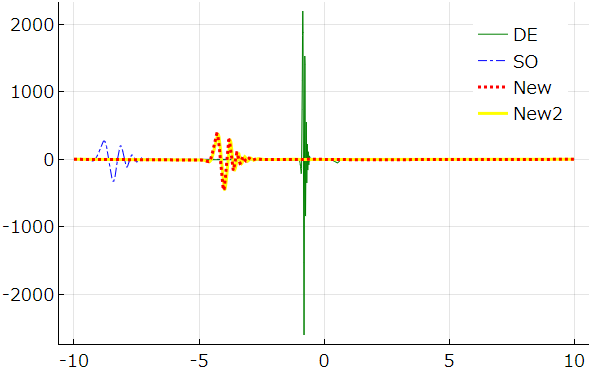}}
	\caption{{\rm (a)} The original integrand $f$.
	{\rm (b)} The transformed integrands $f(\phi(\cdot))\phi '(\cdot)$.
	\label{img:52F}}
\end{figure}

Then, we compare the performances of the formulas $H$
as transformation formulas for integration.
\Cref{table:52} shows the parameters of \Cref{thm:Sugihara_DE},
where the parameter $d$ of DE is calculated as 
\begin{align}
	d_{\mathrm{DE}}
	= \min_{j=1, \ldots, 4}\ \mathrm{Im} \left[
	\sinh^{-1}\!\left(\frac{2}{\pi}\sinh^{-1}
	(\delta_j + \epsilon_j\mathrm{i})\right) \right] .
\end{align}
\Cref{img:52F} shows the original and transformed integrands.
\begin{table}[ht]
	\footnotesize
	\centering
	\caption{Parameters of \Cref{thm:Sugihara_DE}. \label{table:52}}
	\begin{tabular}{|c|cccc|}
		\hline
		& DE & SO & New & New2 \\
		\hline
		$\gamma$ & 1 & 1 & 1 & 1 \\
		$d$ & 0.0976 & $\pi/2$ & $\pi/2$ & ?\\
		$\beta_2$ & 1.57 & $ 5.77\times 10^{-6}$ 
			& $5.12\times 10^{-3}$  & $5.12\times 10^{-3}$\\
		\hline
	\end{tabular}
\end{table}

\begin{figure}
	\subfigure[\label{img:52err}]
	{\includegraphics[scale = 0.39]{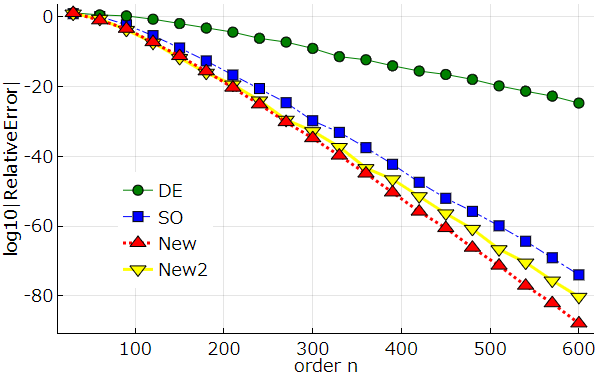}}
	\subfigure[\label{img:52time}]
	{\includegraphics[scale = 0.39]{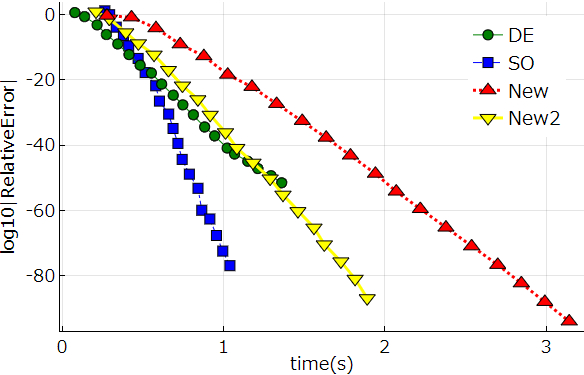}}
	\caption{{\rm (a)} Orders $n$ and errors.
	{\rm (b)} Calculation time and errors.
	The calculation time includes time
	to determine parameters and carry out the trapezoidal formula.}
\end{figure}

Finally, we compare the errors of the numerical integration.
\Cref{img:52err} shows relations between orders $n$ and
the errors. \Cref{img:52time} shows relations between time
for the calculation of numerical integration and the errors.
Here, we assume that $d_{\mathrm{New2}} = \pi/2$
when we calculate the mesh size of the trapezoidal formula \cref{eq:hDE}.

%%%%%%%%%%
\subsection{Integral on a semi-infinite interval (i)}
\label{subsec:53}
%%%%%%%%%%

We consider an integral on a semi-infinite interval:
\begin{align}
	& \int_0^\infty
	\frac
	{\exp\!\left(\frac{1}{50}(\epsilon_1^2 + (x-\delta_1)^2)^{-\frac{3}{2}}\right)\!
	\exp\!\left(\frac{1}{20}(\epsilon_4^2 + (x-\delta_4)^2)^{-\frac{3}{2}}\right)}
	{\sqrt{x}\sqrt{\epsilon_2^2 + (x-\delta_2)^2}\sqrt{\epsilon_3^2 + (x-\delta_3)^2}}
	\mathrm{d}x = 30.6929\ldots ,
\end{align}
where
$\delta_1 \pm \epsilon_1\mathrm{i} = 0.3 \pm 0.2\mathrm{i}$, 
$\delta_2 \pm \epsilon_2\mathrm{i} = 0.5 \pm 0.6\mathrm{i}$, 
$\delta_3 \pm \epsilon_3\mathrm{i} = 0.8 \pm 0.5\mathrm{i}$, 
and $\delta_4 \pm \epsilon_4\mathrm{i} = 1.2 \pm 0.3\mathrm{i}$.
The change of variables for this integral is $\phi (t) = \sinh (H(t))$.
There are singularities as follows:
\begin{align}
	\tilde{S} = \{ \log (\delta_j \pm \epsilon_j \mathrm{i} )
	\mid j = 1, \ldots, 4, \, k \in \mathbb{Z} \}, \quad
	S_\psi = \emptyset .
\end{align}

\begin{figure}[ht]
	\centering
	\subfigure[$H_{\mathrm{SO}}(\mathscr{D}_{\frac{\pi}{2}})$]
	{\includegraphics[scale = 0.26]{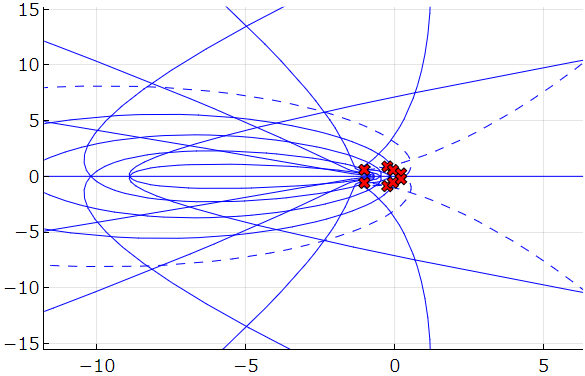}}
	\subfigure[$H_{\mathrm{New}}(\mathscr{D}_{\frac{\pi}{2}})$]
	{\includegraphics[scale = 0.26]{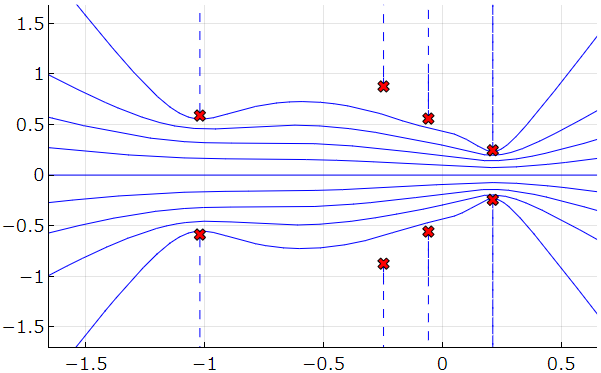}}
	\subfigure[$H_{\mathrm{New2}}(\mathscr{D}_{\frac{\pi}{2}})$]
	{\includegraphics[scale = 0.26]{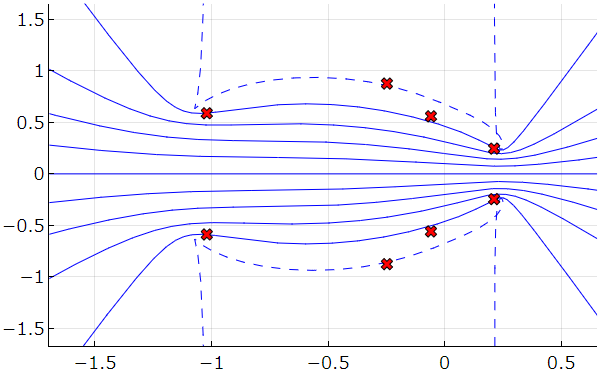}}
	\caption
	{Images $H(\mathscr{D}_{\pi/2})$ and singularities.
	The solid lines are the images of lines
	which are parallel to the real axis in $\mathscr{D}_{\pi/2}$.
	The dotted lines show $H(\partial \mathscr{D}_{\pi/2})$.
	\label{img:53H}}
\end{figure}

First, the formulas $H$ are given by
\begin{align}
	H_{\mathrm{DE}}(t) =&\ \frac{\pi}{2}\sinh(t) , \\
	H_{\mathrm{SO}}(t) \approx&\ 0.784\sinh(t) - 0.894 \\
	& - 1.089\,t - 0.496\,t^2 -0.249\,t^3 \nonumber , \\
	H_{\mathrm{New}}(t) \approx&\ 0.0755\sinh(t+0.693)
	+ 0.492 \tan^{-1}\!\left(\mathrm{e}^{t + 1.91}\right) \\
	& + 0.120 \tan^{-1}\!\left(\mathrm{e}^{t + 1.56}\right)
	+ 0.172 \tan^{-1}\!\left(\mathrm{e}^{t + 0.781}\right)
	- 1.02 \nonumber , \\
	H_{\mathrm{New2}}(t) \approx&\ 0.0755\sinh(t+0.693)
	+ 0.386 \tanh(t+1.91) \\
	& + 0.0944 \tanh(t+1.56)
	+ 0.135 \tanh(t + 0.781)
	- 0.404 \nonumber .
\end{align}
\Cref{img:53H} shows the images $H(\mathscr{D}_{\pi/2})$.

\begin{figure}
	\centering
	\subfigure[]
	{\includegraphics[scale = 0.39]{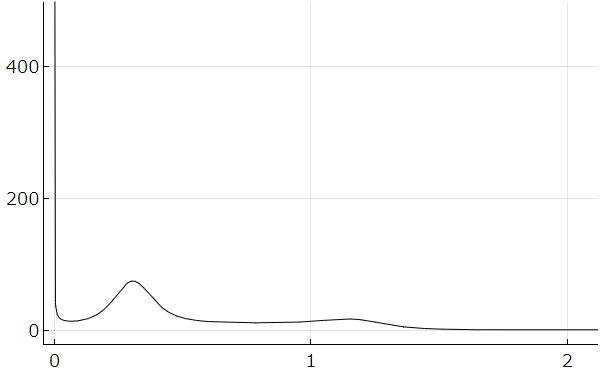}}
	\subfigure[]
	{\includegraphics[scale = 0.39]{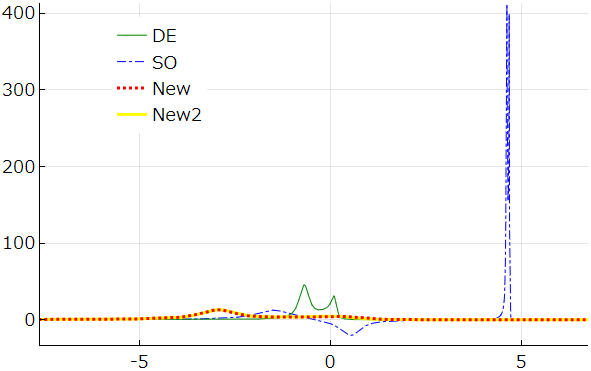}}
	\caption
	{{\rm (a)} The original integrand $f$.
	{\rm (b)} The transformed integrands $f(\phi(\cdot))\phi '(\cdot)$.
	\label{img:53F}}
\end{figure}

Then, we compare the performances of the formulas $H$
as transformation formulas for integration.
\Cref{table:53} shows the parameters of \Cref{thm:Sugihara_DE},
where the parameter $d$ of DE is calculated as 
\begin{align}
	d_{\mathrm{DE}}
	= \min_{j=1, \ldots, 4}\ \mathrm{Im} \left[
	\sinh^{-1}\!\left(\frac{2}{\pi}\log
	(\delta_j + \epsilon_j\mathrm{i})\right) \right] .
\end{align}
\Cref{img:53F} shows the original and transformed integrands.
\begin{table}[ht]
	\footnotesize
	\centering
	\caption{Parameters of \Cref{thm:Sugihara_DE}. \label{table:53}}
	\begin{tabular}{|c|cccc|}
		\hline
		& DE &SO & New & New2 \\
		\hline
		$\gamma$ & 1 & 1 & 1 & 1 \\
		$d$ & 0.155 & ? & $\pi/2$ & ? \\
		$\beta_2$ & 0.393 & 0.196 & 0.0377 & 0.0377 \\
		\hline
	\end{tabular}
\end{table}

\begin{figure}
	\subfigure[\label{img:53err}]
	{\includegraphics[scale = 0.39]{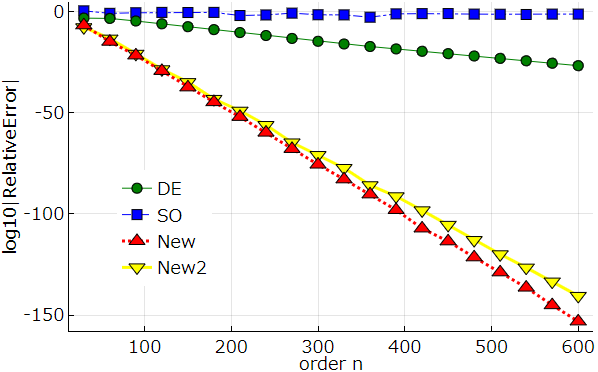}}
	\subfigure[\label{img:53time}]
	{\includegraphics[scale = 0.39]{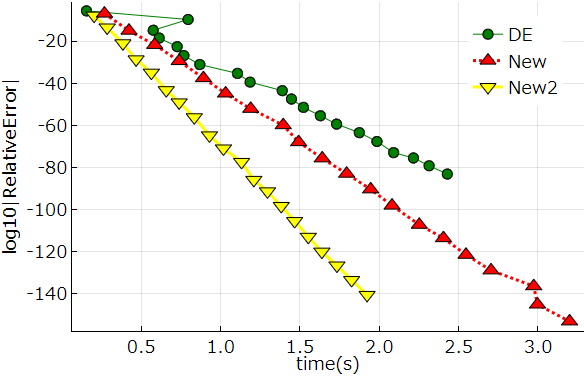}}
	\caption{{\rm (a)} Orders $n$ and errors.
	{\rm (b)} Calculation time and errors.
	The calculation time includes time
	to determine parameters and carry out the trapezoidal formula.
	Results of SO are omitted
	because it took too long to determine the parameters.}
\end{figure}

Finally, we compare the errors of the numerical integration.
\Cref{img:53err} shows relations between orders $n$ and
the errors. \Cref{img:53time} shows relations between time
for the calculation of numerical integration and the errors.
Here, we assume that $d_{\mathrm{SO}}, d_{\mathrm{New2}} = \pi/2$
when we calculate the mesh size of the trapezoidal formula \cref{eq:hDE}.

%%%%%%%%%%
\subsection{Integral on a semi-infinite interval (ii)}
\label{subsec:54}
%%%%%%%%%%

We consider an integral on a semi-infinite interval:
\begin{align}
	& \int_0^\infty \!
	\cos\!\left(\!\frac{5}{\epsilon_1^2 + (x-\delta_1)^2}\!\right) \!
	\cos\!\left(\!\frac{10}{\epsilon_7^2 + (x-\delta_7)^2}\!\right) \!
	\prod_{j=2}^6
	\exp\!\left(\!\frac{c_j}{\epsilon_j^2 + (x-\delta_j)^2}\!\right) \!
	\frac{\mathrm{e}^{-\frac{1}{5}x}}{\sqrt{x}}
	\mathrm{d}x \\
	&= -0.3451\ldots \nonumber ,
\end{align}
where
$\delta_1 \pm \epsilon_1\mathrm{i} = 1 \pm 0.1\mathrm{i}$, 
$\delta_2 \pm \epsilon_2\mathrm{i} = 2 \pm 0.5\mathrm{i}$, 
$\delta_3 \pm \epsilon_3\mathrm{i} = 3 \pm 0.3\mathrm{i}$, 
$\delta_4 \pm \epsilon_4\mathrm{i} = 4 \pm 0.5\mathrm{i}$,
$\delta_5 \pm \epsilon_5\mathrm{i} = 5 \pm 0.2\mathrm{i}$,
$\delta_6 \pm \epsilon_6\mathrm{i} = 6 \pm 0.5\mathrm{i}$,
$\delta_7 \pm \epsilon_7\mathrm{i} = 7 \pm 0.1\mathrm{i}$,
$c_2 = 0.8$, 
$c_3 = 0.2$, 
$c_4 = 0.5$, 
$c_5 = 0.1$, 
and $c_6 = 0.5$.
The change of variables for this integral is $\phi (t) = \log(\exp(H(t))+1)$.
There are singularities as follows:
\begin{align}
	\tilde{S} = \{ \log (\delta_j \pm \epsilon_j \mathrm{i} )
	\mid j = 1, \ldots, 7, \, k \in \mathbb{Z} \}, \quad
	S_\psi = \{  (\pm 1 + 2k) \pi \mathrm{i} \mid k \in \mathbb{Z} \}.
\end{align}

\begin{figure}[ht]
	\centering
	\subfigure[$H_{\mathrm{New}}(\mathscr{D}_{\frac{\pi}{2}})$]
	{\includegraphics[scale = 0.38]{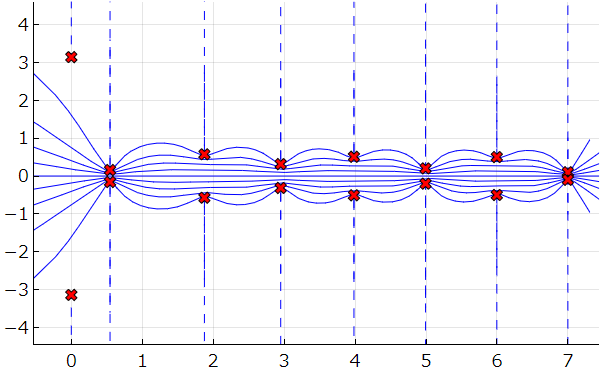}}
	\subfigure[$H_{\mathrm{New2}}(\mathscr{D}_{\frac{\pi}{2}})$]
	{\includegraphics[scale = 0.38]{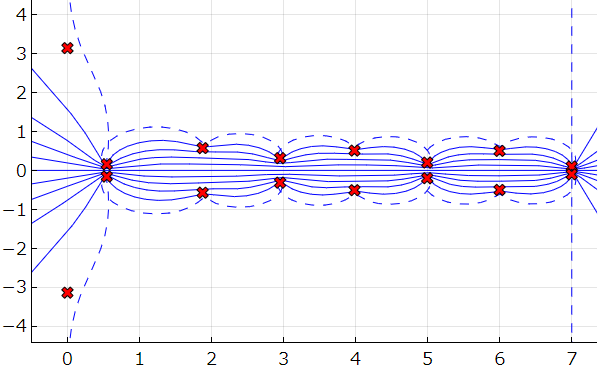}}
	\caption
	{Images $H(\mathscr{D}_{\pi/2})$ and singularities.
	The solid lines are the images of lines
	which are parallel to the real axis in $\mathscr{D}_{\pi/2}$.
	The dotted lines show $H(\partial \mathscr{D}_{\pi/2})$.
	\label{img:54H}}
\end{figure}

First, the formulas $H$ are given by
\begin{align}
	H_{\mathrm{DE}}(t) =& \frac{\pi}{2}\sinh(t), \\
	H_{\mathrm{New}}(t) \approx&\ 1.17\times 10^{-5}\sinh(t+0.458)
	+ 0.348 \tan^{-1}\!\left(\mathrm{e}^{t + 13.4}\right) \\
	& + 0.847 \tan^{-1}\!\left(\mathrm{e}^{t + 7.35}\right)
	+ 0.684 \tan^{-1}\!\left(\mathrm{e}^{t + 5.26}\right) \nonumber \\
	& + 0.657 \tan^{-1}\!\left(\mathrm{e}^{t + 2.08}\right)
	+ 0.642 \tan^{-1}\!\left(\mathrm{e}^{t + 0.0463}\right) \nonumber \\
	& + 0.639 \tan^{-1}\!\left(\mathrm{e}^{t - 3.92}\right)
	+ 0.637 \tan^{-1}\!\left(\mathrm{e}^{t - 5.92}\right), \nonumber \\
	H_{\mathrm{New2}}(t) \approx&\ 1.17\times 10^{-5}\sinh(t+0.458)
	+ 0.273 \tanh(t + 13.4) \\
	& + 0.665 \tanh(t + 7.35) + 0.538 \tanh(t + 5.26) \nonumber \\
	& + 0.516 \tanh(t + 2.08) + 0.504 \tanh(t + 0.0463) \nonumber \\
	& + 0.502 \tanh(t - 3.92) + 0.500 \tanh(t - 5.92) - 3.50 \nonumber,
\end{align}
where the method of Slevinsky and Olver is omitted
because we could not solve the optimization problem
by their program \cite{SOprogram}.
\Cref{img:54H} shows the images $H(\mathscr{D}_{\pi/2})$.

\begin{figure}
	\centering
	\subfigure[]
	{\includegraphics[scale = 0.39]{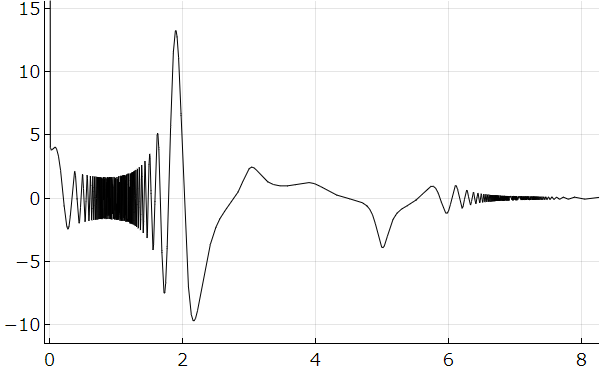}}
	\subfigure[]
	{\includegraphics[scale = 0.39]{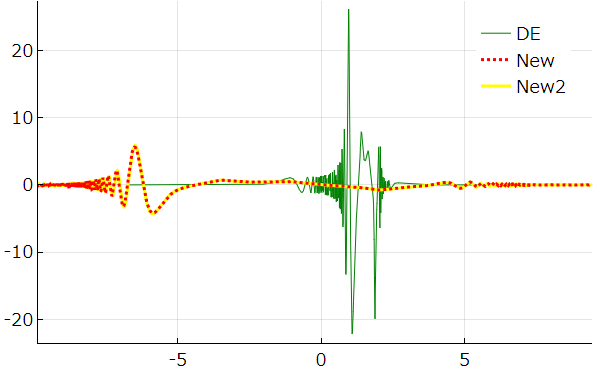}}
	\caption
	{{\rm (a)} The original integrand $f$. 
	{\rm (b)} The transformed integrands $f(\phi(\cdot))\phi '(\cdot)$.
	\label{img:54F}}
\end{figure}

Then, we compare the performances of the formulas $H$
as transformation formulas for integration.
\Cref{table:54} shows the parameters of \Cref{thm:Sugihara_DE},
where the parameter $d$ of DE is calculated as 
\begin{align}
	d_{\mathrm{DE}}
	= \min_{j=1, \ldots, 7}\ \mathrm{Im}\left[
	\sinh^{-1}\left(\frac{2}{\pi}\log(\exp
	(\delta_j + \epsilon_j\mathrm{i}) - 1)\right) \right]
\end{align}
\Cref{img:54F} shows the original and transformed integrands.
\begin{table}[ht]
	\footnotesize
	\centering
	\caption{Parameters of \Cref{thm:Sugihara_DE}. \label{table:54}}
	\begin{tabular}{|c|ccc|}
		\hline
		& DE & New & New2 \\
		\hline
		$\gamma$ & 1 & 1 & 1 \\
		$d$ & 0.0139 & $\pi/2$ & ?\\
		$\beta_2$ & 0.157 & $1.85\times 10^{-6}$ & $1.85\times 10^{-6}$ \\
		\hline
	\end{tabular}
\end{table}

\begin{figure}
	\subfigure[\label{img:54err}]
	{\includegraphics[scale = 0.39]{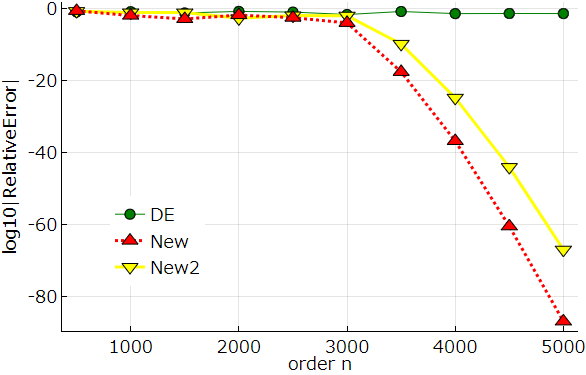}}
	\subfigure[\label{img:54time}]
	{\includegraphics[scale = 0.39]{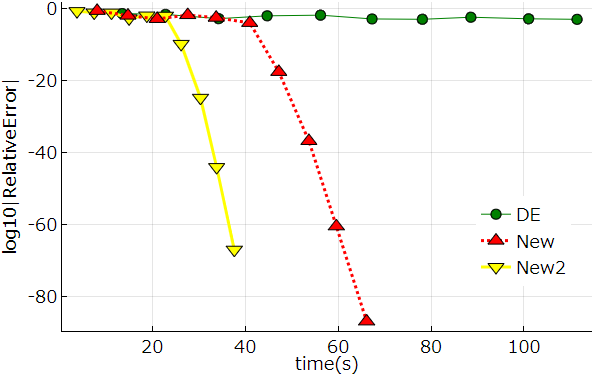}}
	\caption{{\rm (a)} Orders $n$ and errors.
	{\rm (b)} Calculation time and errors.
	The calculation time includes time
	to determine parameters and carry out the trapezoidal formula.}
\end{figure}

Finally, we compare the errors of the numerical integration.
\Cref{img:54err} shows relations between orders $n$ and
the errors. \Cref{img:54time} shows relations between time
for the calculation of numerical integration and the errors.
Here, we assume that $d_{\mathrm{New2}} = \pi/2$
when we calculate the mesh size of the trapezoidal formula \cref{eq:hDE}.

%
%
%%%%%%%%%%%%%%%%%%%%%%%%%%%%%%%%%%%%%%%%%%%%%%%%%%%%%%%%%%%%%%%%%%%%%%%%%%%
\section{Conclusion}
\label{sec:6}
%%%%%%%%%%%%%%%%%%%%%%%%%%%%%%%%%%%%%%%%%%%%%%%%%%%%%%%%%%%%%%%%%%%%%%%%%%%

We improved the DE formula in the case where the integrand has finite singularities.
We improved it by constructing new transformation formulas
$H_{\mathrm{New}}$ and $H_{\mathrm{New2}}$.
The transformation $H_{\mathrm{New}}$ could be considered to be a generalization of
the DE transformations and $H_{\mathrm{New2}}$ was an approximation of it.
By numerical experiments, we confirmed the effectiveness of them
even when the methods of the previous research failed.

In future work, we will need to consider cases
where we do not know the locations of singularities or
integrands have infinite singularities.
\appendix

%%%%%%%%%%%%%%%%%%%%%%%%%%%%%%%%%%%%%%%%%%%%%%%%%%%%%%%%%%%%%%%%%%%%%%%%%%%
\section{Proof of Theorem \ref{thm:asympt}}
\label{ap:asympt}
%%%%%%%%%%%%%%%%%%%%%%%%%%%%%%%%%%%%%%%%%%%%%%%%%%%%%%%%%%%%%%%%%%%%%%%%%%%

We define $I, J, I_0$ and $J_0$ as
\begin{align}
	&I (\alpha_1, \ldots, \alpha_M) = 
	\int_0^t \mathrm{e}^{\Delta\theta \tau}
	\prod_{j=1}^M \cosh^{\alpha_j -1}(\tau-\tau_j)
	\mathrm{d}\tau , \\
	&J (\alpha_1, \ldots, \alpha_M) = \\
	&\int_0^t \mathrm{e}^{\Delta\theta \tau}
	\left\{
		\prod_{j=1}^{M-1} \cosh^{\alpha_j -1}(\tau-\tau_j)
	\right\}
	\cosh^{\alpha_M - 2}(\tau-\tau_M) \sinh(\tau-\tau_M)
	\mathrm{d}\tau , \nonumber \\
	&I_0 =
	\left\{
		\prod_{j=1}^{M-1} \cosh^{\alpha_j -1}(-\tau_j)
	\right\}
	\cosh^{\alpha_M - 2}(-\tau_M)
	\sinh(-\tau_M) , \\
	&J_0 =
	\left\{
		\prod_{j=1}^{M-1} \cosh^{\alpha_j -1}(-\tau_j)
	\right\}
	\cosh^{\alpha_M - 3}(-\tau_M)
	\sinh^2(-\tau_M) .
\end{align}

First, we consider asymptotic expansion of $I$ and $J$. We rearrange them as
\begin{align}
	& I (\alpha_1, \ldots, \alpha_M) \\
	=& \int_0^t \mathrm{e}^{\Delta\theta \tau}
		\left\{
			\prod_{j=1}^{M-1} \cosh^{\alpha_j -1}(\tau-\tau_j)
		\right\}
		\cosh^{\alpha_M - 2}(\tau-\tau_M)
		\left(\sinh(\tau-\tau_M)\strut\right)'
		\mathrm{d}\tau \\
	=&\ \mathrm{e}^{\Delta\theta t}
		\left\{
			\prod_{j=1}^{M-1} \cosh^{\alpha_j -1}(t-\tau_j)
		\right\}
		\cosh^{\alpha_M - 2}(t-\tau_M)
		\sinh(t-\tau_M)  - I_0 \\
	-& \int_0^t
		\Delta\theta\mathrm{e}^{\Delta\theta \tau}
		\left\{
			\prod_{j=1}^{M-1} \cosh^{\alpha_j -1}(\tau-\tau_j)
		\right\}
		\cosh^{\alpha_M - 2}(\tau-\tau_M) \sinh(\tau-\tau_M)
		\mathrm{d}\tau \nonumber \\
	& \begin{aligned}
		-\sum_{j=1}^{M-1}
			&\left[ \vphantom{\prod_{ \substack{k=1 \\ k\neq j} }^{M-1}}\right.
				\int_0^t
				\mathrm{e}^{\Delta\theta \tau}
				(\alpha_j - 1) \cosh^{\alpha_j - 2}(\tau - \tau_j) \sinh(\tau - \tau_j) \\
				&\cdot
				\left\{
					\prod_{ \substack{k=1 \\ k\neq j} }^{M-1}
					\cosh^{\alpha_k - 1}(\tau - \tau_k)
				\right\}
				\cosh^{\alpha_M - 2}(\tau - \tau_M)
				\sinh(\tau - \tau_M)
				\mathrm{d}\tau
				\left. \vphantom{\prod_{ \substack{k=1 \\ k\neq j} }^{M-1}}\right]
	\end{aligned} \nonumber \\
	&\begin{aligned} 
		-\int_0^t
		\mathrm{e}^{\Delta\theta \tau}
		&\left\{
			\prod_{j=1}^{M-1}
			\cosh^{\alpha_j - 1}(\tau-\tau_j)
		\right\} \\
		&\cdot (\alpha_M - 2) \cosh^{\alpha_M - 3}(\tau - \tau_M)
		\sinh^2(\tau - \tau_M)
		\mathrm{d}\tau
	\end{aligned} \nonumber \\
	=&\ \mathrm{e}^{\Delta\theta t} \label{eq:I_expand}
		\left\{
			\prod_{j=1}^{M-1} \cosh^{\alpha_j -1}(t-\tau_j)
		\right\}
		\cosh^{\alpha_M - 2}(t-\tau_M)
		\sinh(t-\tau_M) \\
	&- I_0 -\Delta \theta J (\alpha_1, \ldots, \alpha_M)
	- I(\alpha_1, \ldots, \alpha_M)
	\left\{
		\sum_{j=1}^M (\alpha_j - 1) -1
	\right\} \nonumber \\
	& \begin{aligned}
		+\sum_{j=1}^{M-1}
		[ (\alpha_j - 1) & \cosh(t_j-\tau_M) \\
		& \cdot
		I(\alpha_1, \ldots, \alpha_{j-1}, \alpha_j-1,
		\alpha_{j+1}, \ldots, \alpha_{M-1}, \alpha_M -1) ]
	\end{aligned} \nonumber \\
	&+ 
	(\alpha_M - 2)
	I(\alpha_1, \ldots, \alpha_{M-1}, \alpha_M -2), \nonumber
\end{align}
and
\begin{align}
	& J (\alpha_1, \ldots, \alpha_M) \\
	=& \int_0^t \mathrm{e}^{\Delta\theta \tau}
		\left\{
			\prod_{j=1}^{M-1} \cosh^{\alpha_j -1}(\tau-\tau_j)
		\right\} \\
	&\quad\ \ \cdot \cosh^{\alpha_M - 3}(\tau-\tau_M) \sinh(\tau-\tau_M)
		\left(\sinh(\tau-\tau_M)\strut\right)'
		\mathrm{d}\tau \nonumber \\
	=&\ \mathrm{e}^{\Delta\theta t}
		\left\{
			\prod_{j=1}^{M-1} \cosh^{\alpha_j -1}(t-\tau_j)
		\right\}
		\cosh^{\alpha_M - 3}(t-\tau_M)
		\sinh^2(t-\tau_M) -J_0 \\
	-& \int_0^t \! \!
		\Delta\theta\mathrm{e}^{\Delta\theta \tau} \!
		\left\{ \!
			\prod_{j=1}^{M-1} \cosh^{\alpha_j -1}(\tau-\tau_j)
		\right\} \!
		\cosh^{\alpha_M - 3}(\tau-\tau_M) \sinh^2(\tau-\tau_M)
		\mathrm{d}\tau \nonumber \\
	& \begin{aligned}
		-\sum_{j=1}^{M-1}
			\!\! & \left[ \vphantom{\prod_{ \substack{k=1 \\ k\neq j} }^{M-1}}\right.
			\int_0^t
			\mathrm{e}^{\Delta\theta \tau}
			(\alpha_j - 1) \cosh^{\alpha_j - 2}(\tau-\tau_j) \sinh(\tau-\tau_j) \\
			& \cdot
			\left\{
				\prod_{ \substack{k=1 \\ k\neq j} }^{M-1}
				\cosh^{\alpha_k - 1}(\tau-\tau_k)
			\right\}
			\cosh^{\alpha_M - 3}(\tau-\tau_M)
			\sinh^2(\tau-\tau_M)
			\mathrm{d}\tau
			\left. \vphantom{\prod_{ \substack{k=1 \\ k\neq j} }^{M-1}}\right]
	\end{aligned} \nonumber \\
	&\begin{aligned}
		-\int_0^t
			\mathrm{e}^{\Delta\theta \tau}
			&\left\{ \prod_{j=1}^{M-1} \cosh^{\alpha_j - 1}(\tau-\tau_j)
			\right\} \\
		& \cdot (\alpha_M - 3) \cosh^{\alpha_M - 4}(\tau-\tau_M)
			\sinh^3(\tau-\tau_M)
			\mathrm{d}\tau
	\end{aligned} \nonumber \\
	&- \int_0^t \mathrm{e}^{\Delta\theta \tau}
		\left\{ \prod_{j=1}^{M-1} \cosh^{\alpha_j -1}(\tau-\tau_j)
		\right\}
		\cosh^{\alpha_M - 2}(\tau-\tau_M) \sinh(\tau-\tau_M)
		\mathrm{d}\tau \nonumber \\
	=&\ \mathrm{e}^{\Delta\theta t} \label{eq:J_expand}
		\left\{
			\prod_{j=1}^{M-1} \cosh^{\alpha_j -1}(t-\tau_j)
		\right\}
		\cosh^{\alpha_M - 3}(t-\tau_M)
		\sinh^2(t-\tau_M) -J_0 \\
	&- \Delta\theta
		\left( I(\alpha_1, \ldots, \alpha_M) - I(\alpha_1, \ldots, \alpha_M -2)
		\right)	\nonumber \\
	&- J(\alpha_1, \ldots, \alpha_M)
		\left\{ \sum_{j=1}^M (\alpha_j - 1) -1
		\right\} \nonumber \\
	&\ \begin{aligned} 
		+\sum_{j=1}^{M-1} [
		&(\alpha_j - 1)\cosh(\tau_j-\tau_M) \\
		&\cdot J(\alpha_1, \ldots, \alpha_{j-1}, \alpha_j-1, 
		\alpha_{j+1}, \ldots, \alpha_{M-1}, \alpha_M -1) ]
	\end{aligned} \nonumber \\
	&+ (\alpha_M - 3)
		J(\alpha_1, \ldots, \alpha_{M-1}, \alpha_M -2). \nonumber
\end{align}
Some of the terms in \cref{eq:I_expand} and \cref{eq:J_expand} can be ignored
because they are on the order of $\mathrm{O}(1)$ as $|t| \to \infty$.
Indeed, the polygon $P$ in \Cref{subsec:SC} has $(2M + 2)$ vertices
including $\pm \infty$. Then, we see
\begin{align}
	2\sum_{j=1}^M \alpha_j - (\theta_+ + \theta_-) = 2M
	\quad\Leftrightarrow\quad
	\sum_{j=1}^M (\alpha_j -1) = \bar{\theta},
\end{align}
and specifically,
\begin{align}
	\sum_{j=1}^M (\alpha_j -1) + \Delta\theta = \theta_+, \quad  
	\sum_{j=1}^M (\alpha_j -1) - \Delta\theta = \theta_-.
\end{align}
Since $\theta_+$ and $\theta_-$ satisfy
$0 \leq \theta_+,\, \theta_- \leq 1$, we see that
\begin{multline}
	I (\alpha_1, \ldots, \alpha_{j-1}, \alpha_j-1, 
	\alpha_{j+1},\, \ldots, \alpha_{M-1}, \alpha_M -1) \\
	= \left\{
		\begin{aligned}
			\int_0^t \mathrm{O}\!\left(\mathrm{e}^{(\theta_+ -2)\tau}\right) \mathrm{d}\tau
			&= \mathrm{O}(1) \quad (t \to +\infty) \\
			\int_0^t \mathrm{O}\!\left(\mathrm{e}^{(2 - \theta_-)\tau}\right) \mathrm{d}\tau
			&= \mathrm{O}(1) \quad (t \to -\infty)
		\end{aligned} \right. .
\end{multline}
We can show 
$I (\alpha_1, \ldots, \alpha_{M-1}, \alpha_M -2)$,
$J (\alpha_1, \ldots, \alpha_{j-1}, \alpha_j-1, \alpha_{j+1}, \ldots, \alpha_{M-1}, \alpha_M -1)$, 
and
$J (\alpha_1, \ldots, \alpha_{M-1}, \alpha_M -2) = \mathrm{O}(1) \ (|t| \to \infty)$
in a similar manner.

Furthermore, since the first term of \cref{eq:I_expand} is rearranged as
\begin{multline}
	\mathrm{e}^{\Delta\theta t}
	\left\{
		\prod_{j=1}^{M-1} \cosh^{\alpha_j -1}(t-\tau_j)
	\right\}
	\cosh^{\alpha_M - 2}(t-\tau_M)
	\sinh(t-\tau_M) \\
	= \frac{\mathrm{e}^{\Delta\theta t}}{2^{\bar{\theta}-1}}
	\frac{\mathrm{e}^{\bar{\theta}t-\sum_{j=1}^M (\alpha_j - 1) \tau_j}}{2}
	\left( 1 + \mathrm{O}(\mathrm{e}^{-2t}) \right)
\end{multline}
as $t \to +\infty$ and
\begin{multline}
	\mathrm{e}^{\Delta\theta t}
	\left\{
		\prod_{j=1}^{M-1} \cosh^{\alpha_j -1}(t-\tau_j)
	\right\}
	\cosh^{\alpha_M - 2}(t-\tau_M)
	\sinh(t-\tau_M) \\
	= - \frac{\mathrm{e}^{\Delta\theta t}}{2^{\bar{\theta}-1}}
	\frac{\mathrm{e}^{-\bar{\theta}t+\sum_{j=1}^M (\alpha_j - 1) \tau_j}}{2}
	\left( 1 + \mathrm{O}(\mathrm{e}^{2t}) \right)
\end{multline}
as $t \to -\infty$, we can write
\begin{multline} \label{eq:I_asympt1}
	\mathrm{e}^{\Delta\theta t}
	\left\{
		\prod_{j=1}^{M-1} \cosh^{\alpha_j -1}(t-\tau_j)
	\right\}
	\cosh^{\alpha_M - 2}(t-\tau_M)
	\sinh(t-\tau_M) \\
	= \frac{\mathrm{e}^{\Delta\theta t}}{2^{\bar{\theta}-1}}
	\sinh \! \left( \bar{\theta}t - \sum_{j=1}^M (\alpha_j - 1) \tau_j \right)
	+ \mathrm{O}(1)
\end{multline}
as $|t| \to \infty$. Similarly, the first term of \cref{eq:J_expand} is written as
\begin{multline} \label{eq:J_asympt1}
	\mathrm{e}^{\Delta\theta t}
	\left\{
		\prod_{j=1}^{M-1} \cosh^{\alpha_j -1}(t-\tau_j)
	\right\}
	\cosh^{\alpha_M - 3}(t-\tau_M)
	\sinh^2(t-\tau_M) \\
	= \frac{\mathrm{e}^{\Delta\theta t}}{2^{\bar{\theta}-1}}
	\cosh \left( \bar{\theta}t - \sum_{j=1}^M (\alpha_j - 1) \tau_j \right)
	+ \mathrm{O}(1) \quad (|t| \to \infty,\, t \in \mathbb{R}).
\end{multline}

By solving \cref{eq:I_expand} and \cref{eq:J_expand} with respect to $I$ 
and using \cref{eq:I_asympt1} and \cref{eq:J_asympt1}, we obtain
\begin{align}
	&I (\alpha_1, \ldots, \alpha_M) \\
	=&\ 
	\frac{1}{\theta_+\theta_-}%{\bar{\theta}^2-(\Delta\theta)^2}
	\left[\vphantom{\prod_{j=1}^{M-1}}\right.
		\bar{\theta}
		\mathrm{e}^{\Delta\theta t}
		\left\{
			\prod_{j=1}^{M-1} \cosh^{\alpha_j -1}(t-\tau_j)
		\right\}
		\cosh^{\alpha_M - 2}(t-\tau_M)
		\sinh(t-\tau_M) \\
		&- \Delta\theta
		\mathrm{e}^{\Delta\theta t}
		\left\{
			\prod_{j=1}^{M-1} \cosh^{\alpha_j -1}(t-\tau_j)
		\right\}
		\cosh^{\alpha_M - 3}(t-\tau_M)
		\sinh^2(t-\tau_M)
	\left.\vphantom{\prod_{j=1}^{M-1}}\right] \nonumber \\
	&+ \mathrm{O}(1) \nonumber \\
	=&\ 
	\frac{1}{\theta_+ \theta_-2^{\bar\theta-1}}
	\frac{1}{2}
	\left(
		\theta_-
		\mathrm{e}^{\theta_+ t- \sum_{j=1}^M (\alpha_j - 1) \tau_j}
		- \theta_+
		\mathrm{e}^{\theta_- t+ \sum_{j=1}^M (\alpha_j - 1) \tau_j}
	\right)
	+\mathrm{O}(1).
\end{align}

%%%%%%%%%%%%%%%%%%%%%%%%%%%%%%%%%%%%%%%%%%%%%%%%%%%%%%%%%%%%%%%%%%%%%%%%%%%
\section{Examination with respect to the parameter \textit{C}}
\label{ap:C}
%%%%%%%%%%%%%%%%%%%%%%%%%%%%%%%%%%%%%%%%%%%%%%%%%%%%%%%%%%%%%%%%%%%%%%%%%%%

In this section, we show relations between polygons $P$
and the corresponding parameters $C$ experimentally.
We calculate the parameter $C$ using sc-toolbox, a MATLAB
program to solve the Schwarz-Christoffel parameter problem \cite{sc-toolbox}.

\begin{figure}[ht]
	\centering
	\subfigure[]
	{\includegraphics[scale = 0.21]{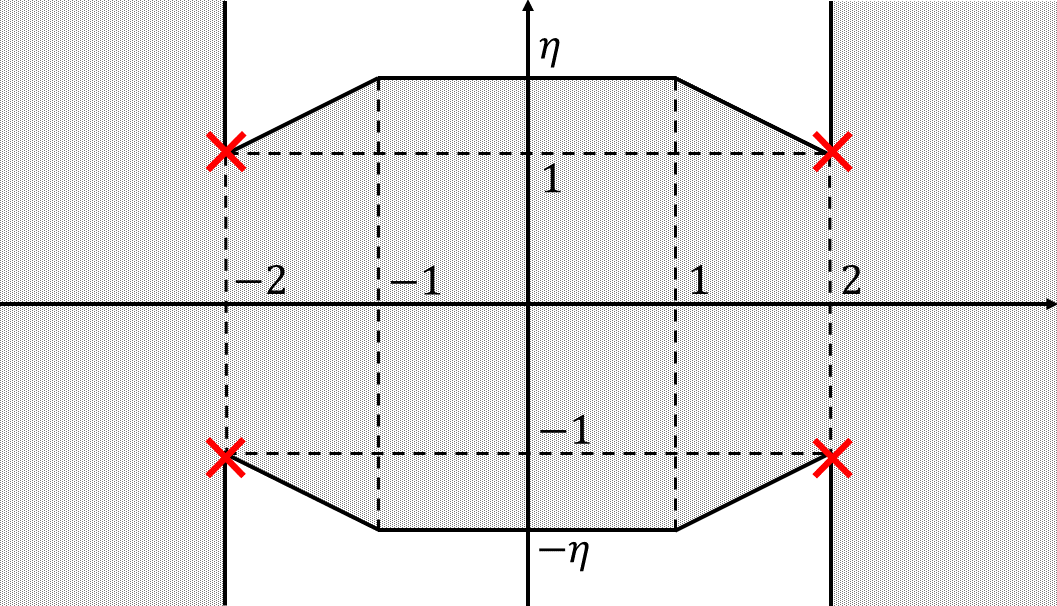}}
	\subfigure[]
	{\includegraphics[scale = 0.21]{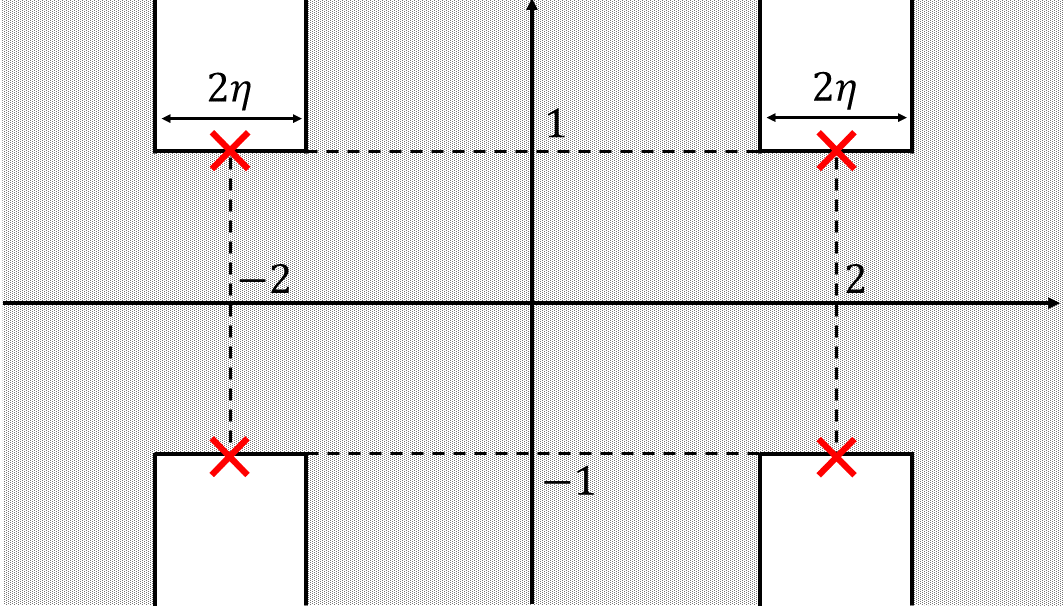}}
	\caption{ Polygons $P$ in the experiment.
	\label{img:polygon_CC}}
\end{figure}

Let $\tilde{\delta}_1 \pm \tilde{\epsilon}_1\mathrm{i} = -2 \pm \mathrm{i}$
and $\tilde{\delta}_2 \pm \tilde{\epsilon}_2\mathrm{i} = 2 \pm \mathrm{i}$
be singularities which the polygon $P$ need to avoid.
Then, let $\eta$ be a positive number and 
we consider the following two kinds of polygons $P$:
\begin{enumerate}
\renewcommand{\labelenumi}{(\alph{enumi})}
	\item
	A polygon which connects the singularities and 4 vertices $(\pm 1, \pm \eta)$.
	\item
	A polygon which avoids the singularities by slits with width $2 \eta$.
\end{enumerate}
We show these polygons in \Cref{img:polygon_CC}.

\begin{figure}[ht]
	\centering
	\subfigure[]
	{\includegraphics[scale = 0.34]{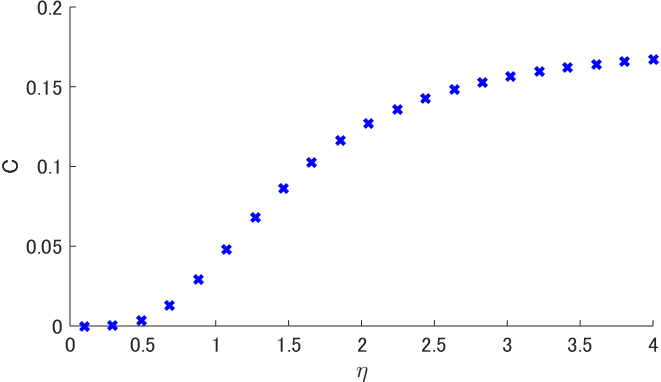}}
	\subfigure[]
	{\includegraphics[scale = 0.34]{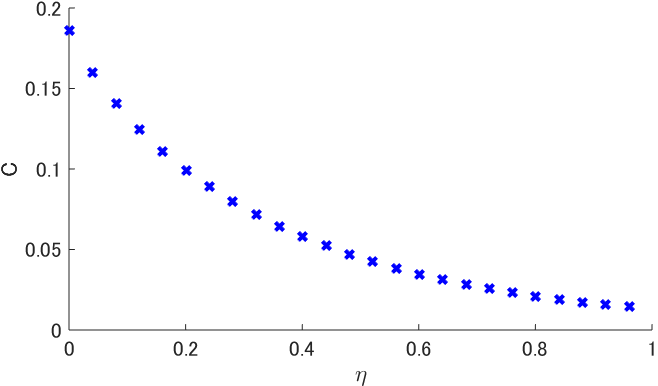}}
	\caption{ Relations between $\eta$ and $C$.
	\label{img:result_CC}}
\end{figure}

\Cref{img:result_CC} shows relations between the parameters $\eta$ and $C$.
We see that the larger the area of $P$ is, the larger the parameter $C$ is.

%%%%%%%%%%%%%%%%%%%%%%%%%%%%%%%%%%%%%%%%%%%%%%%%%%%%%%%%%%%%%%%%%%%%%%%%%%%
\section{Determination of the parameter \textit{T} for the other intervals}
\label{ap:T}
%%%%%%%%%%%%%%%%%%%%%%%%%%%%%%%%%%%%%%%%%%%%%%%%%%%%%%%%%%%%%%%%%%%%%%%%%%%

In the paper, we show how to determine the parameter $T$ in the case
where the interval is $(-1, 1)$.
In this section, we show the other cases.

%%%%%%%%%%
\subsection{Infinite interval}
%%%%%%%%%%

We consider an integral of the interval $(-\infty, \infty)$.
We assume that the integrand $f$ is smooth on the interval $(-\infty, \infty)$
and satisfies
\begin{align}
	f(x) = 
	\left\{
		\begin{aligned}
			& \mathrm{O} \!\left( |x|^r \right) \quad ( x \to +\infty) \\
			& \mathrm{O} \!\left( |x|^s \right) \quad ( x \to -\infty)
		\end{aligned}
	\right.
\end{align}
for some $r, s < -1$.
The change of variables is given as
\begin{align}
	x = \phi(t) = \sinh(H_{\mathrm{New}}(t)).
\end{align}
The decay rate of the transformed integrand is estimated as
\begin{align}
	f(\phi(t))\phi'(t) = 
	\left\{
		\begin{aligned}
			& \mathrm{O} \!\left( \exp\!
				\left( \! \left( \frac{C}{2}(1+r) + \varepsilon \right)
				\mathrm{e}^{t-T} \right) \right)
			\quad (t\to+\infty) \\
			& \mathrm{O} \!\left( \exp\!
				\left( \! \left(\frac{C}{2}(1+s) + \varepsilon \right)
				\mathrm{e}^{T-t} \right) \right)
			\quad (t\to-\infty)
		\end{aligned}
	\right.
\end{align}
for arbitrary $\varepsilon >0$.
Then we see that the parameter $\beta_2$ satisfies
\begin{align}
	\beta_2 \leq
	\min\!
	\left\{
		\left(-\frac{C}{2}(1+r) -\varepsilon \right)\mathrm{e}^{-T},\,
		\left(-\frac{C}{2}(1+s) -\varepsilon \right)\mathrm{e}^{T}
	\right\}.
\end{align}

To make the parameter $\beta_2$ larger,
we we make $\varepsilon$ go to 0 and determine $T$ as
\begin{align}
	-\frac{C}{2}(1+r)\mathrm{e}^{-T} = -\frac{C}{2}(1+s)\mathrm{e}^T
	\quad \Leftrightarrow \quad 
	T = \frac{1}{2}
	\log \left( \frac{1+r}{1+s} \right).
\end{align}
Then the supremum of the parameter $\beta_2$ is estimated as
\begin{align}
	\beta_2^* = \frac{C}{2}\sqrt{(r+1)(s+1)}.
\end{align}

%%%%%%%%%%%
\subsection{Semi-infinite interval (i)}
%%%%%%%%%%%

We consider an integral of the interval $(0, \infty)$.
We assume that the integrand $f$ is smooth on the interval $(0, \infty)$
and satisfies
\begin{align}
	f(x) = 
	\left\{
		\begin{aligned}
			& \mathrm{O} \!\left( x^r \right) \quad ( x \to +\infty) \\
			& \mathrm{O} \!\left( x^q \right) \quad ( x \to +0)
		\end{aligned}
	\right.
\end{align}
for some $r < -1$ and $q > -1$.
The change of variables is given as
\begin{align}
	x = \phi(t) = \exp(H_{\mathrm{New}}(t)).
\end{align}
The decay rate of the transformed integrand is estimated as
\begin{align}
	f(\phi(t))\phi'(t) = 
	\left\{
		\begin{aligned}
			& \mathrm{O} \!\left( \exp\!
				\left( \! \left( \frac{C}{2}(1+r) + \varepsilon \right)
				\mathrm{e}^{t-T} \right) \right)
				\quad (t\to+\infty) \\
			& \mathrm{O} \!\left( \exp\!
				\left( \! - \left(\frac{C}{2}(1+q) - \varepsilon \right)
				\mathrm{e}^{T-t} \right) \right)
				\quad (t\to-\infty)
		\end{aligned}
	\right.
\end{align}
for arbitrary $\varepsilon >0$.
Then we see that the parameter $\beta_2$ satisfies
\begin{align}
	\beta_2 \leq
	\min\!
	\left\{
		\left(-\frac{C}{2}(1+r) -\varepsilon \right)\mathrm{e}^{-T},\,  
		\left(\frac{C}{2}(1+q) -\varepsilon \right)\mathrm{e}^{T}
	\right\}.
\end{align}

To make the parameter $\beta_2$ larger,
we we make $\varepsilon$ go to 0 and determine $T$ as
\begin{align}
	-\frac{C}{2}(1+r)\mathrm{e}^{-T} = \frac{C}{2}(1+q)\mathrm{e}^T
	\quad \Leftrightarrow \quad 
	T = \frac{1}{2}
	\log\! \left(\! - \frac{1+r}{1+q} \right).
\end{align}
Then the supremum of the parameter $\beta_2$ is estimated as
\begin{align}
	\beta_2^* = \frac{C}{2}\sqrt{-(1+r)(1+q)}.
\end{align}

%%%%%%%%%%
\subsection{Semi-infinite interval (ii)}
%%%%%%%%%%

We consider an integral of the interval $(0, \infty)$.
We assume that the integrand $f$ is smooth on the interval $(0, \infty)$
and satisfies
\begin{align}
	f(x) = 
	\left\{
		\begin{aligned}
			& \mathrm{O} \!\left( \mathrm{e}^{-(v-\varepsilon)x} \right)
			\quad ( x \to +\infty) \\
			& \mathrm{O} \!\left( x^q \right) \quad ( x \to +0)
		\end{aligned}
	\right.
\end{align}
for some $q > -1, v  > 0$, and arbitral $\varepsilon > 0$.
The change of variables is given as
\begin{align}
	x = \phi(t) = \log (\exp (H_{\mathrm{New}}(t))+1).
\end{align}
The decay rate of the transformed integrand is estimated as
\begin{align}
	f(\phi(t))\phi'(t) =
	\left\{
		\begin{aligned}
			& \mathrm{O}\!
			\left( \exp \!
				\left(-\frac{C}{2}(v-\varepsilon)\mathrm{e}^{t-T} \right)
			\right)
			\quad (t\to+\infty) \\
			& \mathrm{O}\!
			\left( \exp\!
				\left(-\frac{C}{2}((1+q)-\varepsilon)\mathrm{e}^{T-t}\right)
			\right)
			\quad (t\to-\infty)
		\end{aligned}
	\right.
\end{align}
for arbitrary $\varepsilon >0$.
Then we see that the parameter $\beta_2$ satisfies
\begin{align}
	\beta_2 \leq
	\min\!
	\left\{
		\frac{C}{2}(v-\varepsilon)\mathrm{e}^{-T},\,  
		\frac{C}{2}((1+q)-\varepsilon)\mathrm{e}^{T}
	\right\}.
\end{align}

To make the parameter $\beta_2$ larger,
we we make $\varepsilon$ go to 0 and determine $T$ as
\begin{align}
	\frac{C}{2}v\mathrm{e}^{-T} = \frac{C}{2}(1+q)\mathrm{e}^T
	\quad \Leftrightarrow \quad 
	T = \frac{1}{2}
	\log\! \left( \frac{v}{1+q} \right).
\end{align}
Then the supremum of the parameter $\beta_2$ is estimated as
\begin{align}
	\beta_2^* = \frac{C}{2}\sqrt{v(1+q)}.
\end{align}

%%%%%%%%%%%%%%%%%%%%%%%%%%%%%%%%%%%%%%%%%%%%%%%%%%%%%%%%%%%%%%%%%%%%%%%%%%%
\section{Proof of Theorem \ref{thm:ababa}}
\label{ap:ababa}
%%%%%%%%%%%%%%%%%%%%%%%%%%%%%%%%%%%%%%%%%%%%%%%%%%%%%%%%%%%%%%%%%%%%%%%%%%%

For simplicity, we show the case of $m=4$.

We rearrange the given system of equations as follows:
\begin{align}
		&\left\{
			\begin{aligned}	
				C\sinh(a_1-T) - \frac{D_1}{\sinh(a_1-b_1)}
				- \frac{D_2}{\sinh(a_1-b_2)} - \frac{D_3}{\sinh(a_1-b_3)} = 0 \\
				C\sinh(a_2-T) - \frac{D_1}{\sinh(a_2-b_1)}
				- \frac{D_2}{\sinh(a_2-b_2)} - \frac{D_3}{\sinh(a_2-b_3)} = 0 \\
				C\sinh(a_3-T) - \frac{D_1}{\sinh(a_3-b_1)}
				- \frac{D_2}{\sinh(a_3-b_2)} - \frac{D_3}{\sinh(a_3-b_3)} = 0 \\
				C\sinh(a_4-T) - \frac{D_1}{\sinh(a_4-b_1)}
				- \frac{D_2}{\sinh(a_4-b_2)} - \frac{D_3}{\sinh(a_4-b_3)} = 0
			\end{aligned}
		\right. \\
		\Leftrightarrow \nonumber \\
		\label{eq:Hequations_dash}
		&\left\{
			\begin{aligned}
				C' & (A_1 - T')(A_1 - B_1)(A_1 - B_2)(A_1 - B_3)
				- D_1' A_1(A_1 - B_2)(A_1 - B_3) \\
				& - D_2' A_1(A_1 - B_1)(A_1 - B_3)
				- D_3' A_1(A_1 - B_1)(A_1 - B_2) = 0 \\
				C' & (A_2 - T')(A_2 - B_1)(A_2 - B_2)(A_2 - B_3)
				- D_1' A_2(A_2 - B_2)(A_2 - B_3) \\
				& - D_2' A_2(A_2 - B_1)(A_2 - B_3)
				- D_3' A_2(A_2 - B_1)(A_2 - B_2) = 0 \\
				C' & (A_3 - T')(A_3 - B_1)(A_3 - B_2)(A_3 - B_3)
				- D_1' A_3(A_3 - B_2)(A_3 - B_3) \\
				& - D_2' A_3(A_3 - B_1)(A_3 - B_3)
				- D_3' A_3(A_3 - B_1)(A_3 - B_2) = 0 \\
				C' & (A_4 - T')(A_4 - B_1)(A_4 - B_2)(A_4 - B_3)
				- D_1' A_4(A_4 - B_2)(A_4 - B_3) \\
				& - D_2' A_4(A_4 - B_1)(A_4 - B_3)
				- D_3' A_4(A_4 - B_1)(A_4 - B_2) = 0,
			\end{aligned}
		\right.
\end{align}
where 
$A_1 = \mathrm{e}^{2a_1}, A_2 = \mathrm{e}^{2a_2}, 
A_3 = \mathrm{e}^{2a_3}, A_4 = \mathrm{e}^{2a_4},
B_1 = \mathrm{e}^{2b_1}, B_2 = \mathrm{e}^{2b_2}, 
B_3 = \mathrm{e}^{2b_3}, T' = \mathrm{e}^{2T},  
C' = \frac{1}{16}\mathrm{e}^{-T-b_1-b_2-b_3}, D'_1 = \frac{1}{4}\mathrm{e}^{-b_2-b_3}, 
D'_2 = \frac{1}{4}\mathrm{e}^{-b_1-b_3}$, and  $D'_3 = \frac{1}{4}\mathrm{e}^{-b_1-b_2}$.
The equations \cref{eq:Hequations_dash} can be seen as a linear system of equations
with respect to $C', D'_1, D'_2$, and $D'_3$.
Since it has non-trivial solutions, we see that $\det X = 0$, where
\begin{multline}
	X = \left[
		\begin{array}{cc}
			(A_1 - T')(A_1 - B_1)(A_1 - B_2)(A_1 - B_3) & A_1(A_1 - B_2)(A_1 - B_3) \\
			(A_2 - T')(A_2 - B_1)(A_2 - B_2)(A_1 - B_3) & A_2(A_2 - B_2)(A_2 - B_3) \\
			(A_3 - T')(A_3 - B_1)(A_3 - B_2)(A_1 - B_3) & A_3(A_3 - B_2)(A_3 - B_3) \\
			(A_4 - T')(A_4 - B_1)(A_4 - B_2)(A_1 - B_3) & A_4(A_4 - B_2)(A_4 - B_3)
		\end{array}
	\right. \\
	\left.
		\begin{array}{cc}
			A_1(A_1 - B_1)(A_1 - B_3) & A_1(A_1 - B_1)(A_1 - B_2) \\
			A_2(A_2 - B_1)(A_2 - B_3) & A_2(A_2 - B_1)(A_2 - B_2) \\
			A_3(A_3 - B_1)(A_3 - B_3) & A_3(A_3 - B_1)(A_3 - B_2) \\
			A_4(A_4 - B_1)(A_4 - B_3) & A_4(A_4 - B_1)(A_4 - B_2)
		\end{array}
	\right] .
\end{multline}

Here, we define $X_0$ as
\begin{multline}
	X_0 = \left[ 
		\begin{array}{cc}
			(A_1 - B_1)(A_1 - B_2)(A_1 - B_3) & (A_1 - B_2)(A_1 - B_3) \\
			(A_2 - B_1)(A_2 - B_2)(A_2 - B_3) & (A_2 - B_2)(A_2 - B_3) \\
			(A_3 - B_1)(A_3 - B_2)(A_3 - B_3) & (A_3 - B_2)(A_3 - B_3) \\
			(A_4 - B_1)(A_4 - B_2)(A_4 - B_3) & (A_4 - B_2)(A_4 - B_3)
		\end{array}
	\right. \\
	\left.
		\begin{array}{cc}
			(A_1 - B_1)(A_1 - B_3) & (A_1 - B_1)(A_1 - B_2) \\
			(A_2 - B_1)(A_2 - B_3) & (A_2 - B_1)(A_2 - B_2) \\
			(A_3 - B_1)(A_3 - B_3) & (A_3 - B_1)(A_3 - B_2) \\
			(A_4 - B_1)(A_4 - B_3) & (A_4 - B_1)(A_4 - B_2)
		\end{array}
	\right].
\end{multline}
Then, from the properties of the determinant, we see that
\begin{align} \label{eq:XandX0}
	\det X = ( A_1 A_2 A_3 A_4 - B_1 B_2 B_3 T' ) \det X_0 = 0.
\end{align}
Also, since $\det X_0$ is a polynomial of degree 9 $( = (m-1)^2)$
and is divisible by $(A_i - A_j)$ and $(B_i - B_j)$ for arbitral $i < j$,
we can write
\begin{align}
	\det X_0 = x_0
	\prod_{\substack{i,j = 1,2,3,4 \\ i <j}}(A_i - A_j)
	\prod_{\substack{i,j = 1,2,3 \\ i <j}}(B_i - B_j)
\end{align}
for some real number $x_0$.

We show that $x_0 \neq 0$ as follows.
Letting $B_1 = A_2$, $B_2 = A_3$, and $B_3 = A_4$ formally,
the matrix $X_0$ is an upper triangular matrix of which
the diagonal components are 
$(A_1-A_2)(A_1-A_3)(A_1-A_4), (A_2-A_3)(A_2-A_4), (A_3-A_2)(A_3-A_4)$,
and $(A_4-A_2)(A_4-A_3)$.
From this reason, we see that 
$\det X_0 \not\equiv 0$, which implies $x_0 \neq 0$.

Therefore, from \cref{eq:XandX0}, we see that 
\begin{align}
	A_1 A_2 A_3 A_4 - B_1 B_2 B_3 T' = 0,
\end{align}
which implies that $a_1 - b_1 + \dots - b_3 + a_4 = T$.

\bibliographystyle{plain}
\bibliography{geDEreferences_arXiv}

\begin{thebibliography}{10}

\bibitem{sc-toolbox}
Toby Driscoll, Everett Kropf, and Alfa Heryudono.
\newblock sc-toolbox.
\newblock https:\slash\slash github.com\slash tobydriscoll\slash sc-toolbox.

\bibitem{modSC}
Louis~H. Howell and Lloyd~N. Trefethen.
\newblock A modified {Schwarz}-{Christoffel} transformation for elongated
  regions.
\newblock {\em SIAM J. Sci. Stat. Comput.}, 11:928--949, 1990.

\bibitem{geDEprogram}
Shunki Kyoya.
\newblock Generalized{D}{E}.
\newblock https://github.com/ShunkiKyoya/generalizedDE.

\bibitem{NLsolve}
Patrick~Kofod Mogensen and Kristoffer Carlson.
\newblock Nlsolve.
\newblock https:\slash\slash github.com\slash JuliaNLSolvers\slash NLsolve.jl.

\bibitem{DEindifinite}
Mayinur Muhammad and Masatake Mori.
\newblock Double exponential formulas for numerical indefinite integration.
\newblock {\em Journal of Conputational and Applied Mathematics}, 161:431--448,
  2003.

\bibitem{SOprogram}
Richard~Mikael Slevinsky and Phil Gaudreau.
\newblock Dequadrature.
\newblock https:\slash\slash github.com\slash MikaelSlevinsky\slash
  DEQuadrature.jl.

\bibitem{SO}
Richard~Mikael Slevinsky and Sheehan Olver.
\newblock On the use of conformal maps for the acceleration convergence of the
  trapezoidal rule and sinc numerical methods.
\newblock {\em SIAM J. Sci. Comput.}, 37(2):A676--A700, 2015.

\bibitem{Sugihara}
Masaaki Sugihara.
\newblock Optimality of the double exponential formula –functional analysis
  approach–.
\newblock {\em Numer. Math.}, 75:379--395, 1997.

\bibitem{DE}
Hidetosi Takahasi and Masatake Mori.
\newblock Double exponential formulas for numerical integration.
\newblock {\em Publ. RIMS. Kyoto Univ.}, 9:721--741, 1974.

\bibitem{Trefethen}
Lloyd~Nicholas Trefethen.
\newblock Numerical computation of the {Schwarz}-{Christoffel} transformation.
\newblock {\em SIAM J. Sci. Stat. Comput.}, 1:82--102, 1980.

\end{thebibliography}
\end{document}